\documentclass[a4paper,12pt]{article}
\usepackage{amssymb,fullpage}

\newtheorem{thm}{Theorem}[section]

\newtheorem{cor}[thm]{Corollary}

\newtheorem{lem}[thm]{Lemma}
\newtheorem{propn}[thm]{Proposition}

\newenvironment{proof}{\\{\bf Proof:}}{\hspace*{\fill}$\square$\endlist\par\bigskip}

\begin{document}

\title{Volume growth and heat kernel estimates\\for the continuum random tree}
\author{David Croydon\footnote{Dept of Statistics, University of Warwick, Coventry, CV4 7AL, UK; \underline{d.a.croydon@warwick.ac.uk.}}\\ \tiny{UNIVERSITY OF WARWICK}}
\maketitle

\begin{abstract}
In this article, we prove global and local (point-wise) volume and heat kernel bounds for the continuum random tree. We demonstrate that there are almost-surely logarithmic global fluctuations and log-logarithmic local fluctuations in the volume of balls of radius $r$ about the leading order polynomial term as $r\rightarrow0$. We also show that the on-diagonal part of the heat kernel exhibits corresponding global and local fluctuations as $t\rightarrow0$ almost-surely. Finally, we prove that this quenched (almost-sure) behaviour contrasts with the local annealed (averaged over all realisations of the tree) volume and heat kernel behaviour, which is smooth.
\end{abstract}

\section{Introduction}
\setcounter{thm}{0}

The continuum random tree (CRT) has, since its introduction by Aldous in \cite{Aldous1}, become an important object in modern probability theory. As well as being the scaling limit of a variety of discrete tree-like objects, see \cite{Aldous1}, \cite{Aldous3}, by a suitable random embedding into $\mathbb{R}^d$, it is possible to describe the support of the integrated super-Brownian excursion (ISE) using the CRT (\cite{Aldous4}). In this article, we investigate and present various bounds for the fundamental solution of the heat equation on the CRT, which can of course be thought of as the transition density of the Brownian motion on the CRT. In the course of doing so, we also prove a collection of asymptotic results about the volume of balls in the CRT. With growing evidence (\cite{HaraSlade}) to support the fact that in high dimensions the incipient infinite cluster of percolation at criticality scales to the ISE, we hope that these results will eventually contribute to the understanding of the asymptotic behaviour of random walks on these lattice objects.

We shall denote by $\mathcal{T}$ the CRT, which is a random set defined on an underlying probability space with probability measure $\mathbf{P}$ (we shall write $\mathbf{E}$ for the expectation under $\mathbf{P}$). It has a natural metric, $d_\mathcal{T}$, and a natural volume measure, $\mu$. The existence of Brownian motion on $\mathcal{T}$, as defined in Section 5.2 of \cite{Aldous2}, has already been proved by Krebs, who constructed a process via a Dirichlet form on $\mathcal{T}$, which was defined as the limit of differential operators, and then Brownian motion was a time change of this process, see \cite{Krebs}. We provide an alternative construction, using the resistance form techniques developed by Kigami in \cite{Kigamidendrite} to define a local, regular Dirichlet form on the measure-metric space $(\mathcal{T},d_\mathcal{T},\mu)$. Given this Dirichlet form, standard results allow the construction of an associated Markov process, $X=(X_t)_{t\geq 0}$ with invariant measure $\mu$, which we show is actually Brownian motion (as defined in \cite{Aldous2}) on $\mathcal{T}$. The construction used here seems more natural, allowing us to define the Dirichlet form for Brownian motion directly. Furthermore, the arguments we use to deduce our process satisfies the properties of Brownian motion are more concise, using more recently developed techniques for resistance forms, rather than limiting arguments.

Once a Markov process is established on $\mathcal{T}$, it is natural to ask whether it has a transition density, and if it does, what form does the transition density take? The current literature on measure-metric spaces equipped with a resistance form seems to indicate that an important part of the answer to this question is the volume growth of the space with respect to the resistance metric, see \cite{BarlowBassKumagai} and \cite{Kumagai}. However, for certain random subsets it has been shown that they do not satisfy the kind of uniform volume growth that is often assumed. For example, in \cite{HamJon}, Hambly and Jones show how for a certain class of random recursive fractals the volume of balls of radius $r$ have fluctuations of order of powers of $\ln r^{-1}$ about a leading order polynomial term, $r^\alpha$. With the random self-similarity of the CRT (\cite{Aldous5}), it is reasonable to expect similar behaviour for the CRT. Indeed, it has already been shown that the CRT and a class of recursive fractals do exhibit the same form of Hausdorff measure function, $r^{\alpha}(\ln\ln r^{-1})^{\theta}$, with $\alpha=2$, $\theta=1$ in the case of the CRT (see \cite{treemeas}, Corollary 1.2 and \cite{GMW}, Theorem 5.2). Motivated by the random fractal examples, heat kernel bounds for measure-metric spaces with volume fluctuations have been established in \cite{Croydon}; these suggest that to establish global heat kernel estimates for the CRT, it will be useful to determine the volume growth behaviour.

Henceforth, define the ball of radius $r$ around the point $\sigma\in\mathcal{T}$ to be
$$B(\sigma,r):=\{\sigma':\:d_\mathcal{T}(\sigma,\sigma')<r\}.$$
In the annealed case (Theorem \ref{annealedstate}), we calculate the volume of a ball of radius $r$ around the root exactly. The expression we obtain is easily seen to be asymptotically equal to $2r^2$ as $r\rightarrow 0$. In the quenched case (Theorem \ref{globalstate}), the behaviour is not as smooth and we see fluctuations in the volume growth of logarithmic order, which confirm the expectations of the previous paragraph. Although it is tight enough to demonstrate the order of the fluctuations, we remark that the upper bound for $\inf_{\sigma\in\mathcal{T}}\mu(B(\sigma,r))$ is almost certainly not optimal (as a consequence, neither are the corresponding lower heat kernel bounds). We conjecture that, up to constants, the lower bound for this quantity is sharp.

{\thm \label{annealedstate} Let $\rho$ be the root of $\mathcal{T}$, then
$$\mathbf{E}(\mu(B(\rho,r)))=1-e^{-2r^2},\hspace{20pt}\forall r\geq 0.$$}

\vspace{-5pt}

{\thm \label{globalstate} $\mathbf{P}$-a.s., there exist random constants $C_1$, $C_2$, $C_3$, $C_4\in (0,\infty)$ such that
$$C_1 r^2 \ln_1 r^{-1} \leq \sup_{\sigma\in\mathcal{T}}\mu (B(\sigma, r))\leq C_2 r^2 \ln_1 r^{-1},$$
and
$$C_3 r^2 \left(\ln_1 r^{-1}\right)^{-1} \leq \inf_{\sigma\in\mathcal{T}}\mu (B(\sigma, r))\leq C_4 r^2 \ln_1\ln_1 r^{-1},$$
for $r\in(0,\mathrm{diam}\mathcal{T})$, where $\mathrm{diam}\mathcal{T}$ is the diameter of $(\mathcal{T},d_{\mathcal{T}})$ and $\ln_1x:=\ln x\vee 1$.}

\hspace{10pt}

Locally, we prove the following volume asymptotics, which show that the volume growth of a ball around a particular point has fluctuations about $r^2$ of the order of $\ln \ln r^{-1}$ asymptotically. This exactly mirrors the $\ln\ln r^{-1}$ local fluctuations exhibited by the random recursive fractals of \cite{HamJon}. We remark that the $\limsup$ result has also been proved (up to a constant multiple) in the course of deriving the Hausdorff measure function of $\mathcal{T}$ in \cite{treemeas}. However, we extend the result proved there by deducing the exact value of the constant.
{\thm \label{localstate} $\mathbf{P}$-a.s., we have
$$\limsup_{r\rightarrow 0}\frac{\mu(B(\sigma,r))}{r^2 \ln\ln r^{-1}}=\frac{8}{\pi^2},$$
and also
$$\liminf_{r\rightarrow 0}\frac{\mu(B(\sigma,r))}{r^2 \left(\ln\ln r^{-1}\right)^{-1}}=2,$$
for $\mu$-a.e. $\sigma\in\mathcal{T}$.}

\hspace{10pt}

The global volume bounds of Theorem \ref{globalstate} mean that the CRT satisfies the non-uniform volume doubling of \cite{Croydon}, $\mathbf{P}$-a.s. Results of that article immediately allow us to deduce the existence of a transition density for the Brownian motion on $\mathcal{T}$ and the following bounds upon it.

{\thm \label{heatkernel} $\mathbf{P}$-a.s., the Brownian motion $X=(X_t)_{t\geq 0}$ on $\mathcal{T}$ exists, and furthermore, it has a transition density $(p_t(\sigma,\sigma'))_{\sigma,\sigma'\in\mathcal{T},t>0}$ that satisfies, for some random constants $C_5,C_6,C_7,C_8$, $t_0>0$ and deterministic $\theta_1,\theta_2,\theta_3\in(0,\infty)$,
\begin{equation}\label{lowerbound1}
p_t(\sigma,\sigma')\geq C_5 t^{-\frac{2}{3}}(\ln_1 t^{-1})^{-\theta_1}\exp \left\{-C_6 \left(\frac{d^{3}}{t}\right)^{1/2}\ln_1 \left(\frac{d}{t}\right)^{\theta_2}\right\},
\end{equation}
and
\begin{equation}\label{upperbound1}
p_t(\sigma,\sigma')\leq C_7 t^{-\frac{2}{3}}(\ln_1 t^{-1})^{1/3}\exp \left\{-C_8 \left(\frac{d^{3}}{t}\right)^{1/2}\ln_1 \left(\frac{d}{t}\right)^{-\theta_3}\right\},
\end{equation}
for all $\sigma,\sigma'\in\mathcal{T}$, $t\in(0,t_0)$, where $d:=d_{\mathcal{T}}(\sigma,\sigma')$ and $\ln_1 x := \ln x\vee 1$.}

\hspace{10pt}

This result demonstrates that the heat kernel decays exponentially away from the diagonal and there can be spatial fluctuations of no more than logarithmic order. The following theorem that we prove for the on-diagonal part of the heat kernel shows that global fluctuations of this order do actually occur. Locally, the results we obtain are not precise enough to demonstrate the $\mathbf{P}$-a.s. existence of fluctuations, see Theorem \ref{localheatkernelstate}. However, they do show that there can only be fluctuations of log-logarithmic order, and combined with the annealed result of Proposition \ref{annealedheatkernelstate}, they prove that log-logarithmic fluctuations occur with positive probability.

{\thm \label{globalheatkernelstate} $\mathbf{P}$-a.s., there exist random constants $C_9, C_{10}, C_{11}, C_{12}, t_1>0$ and deterministic $\theta_4\in(0,\infty)$ such that for all $t\in(0,t_1)$,
\begin{equation}\label{globalheatkernelsup}
C_9 t^{-2/3} (\ln \ln t^{-1})^{-14}\leq\sup_{\sigma\in\mathcal{T}}p_t(\sigma,\sigma)\leq C_{10} t^{-2/3} (\ln t^{-1})^{1/3},
\end{equation}
\begin{equation}\label{globalheatkernelinf}
C_{11} t^{-2/3} (\ln t^{-1})^{-\theta_4}\leq \inf_{\sigma\in\mathcal{T}}p_t(\sigma,\sigma)\leq C_{12} t^{-2/3} (\ln t^{-1})^{-1/3}.
\end{equation}}

{\thm \label{localheatkernelstate} $\mathbf{P}$-a.s., for $\mu$-a.e. $\sigma\in\mathcal{T}$, there exist random constants $C_{13}, C_{14}$, $t_2>0$ such that for all $t\in(0,t_2)$,
$$C_{13} t^{-2/3} (\ln \ln t^{-1})^{-14}\leq p_t(\sigma,\sigma) \leq C_{14} t^{-2/3} (\ln\ln t^{-1})^{1/3},$$
and also
$$\liminf_{t\rightarrow 0}\frac {p_t(\sigma,\sigma)}{t^{-2/3} (\ln\ln t^{-1})^{-1/3}}<\infty.$$}

The final estimates we prove are annealed heat kernel bounds at the root of $\mathcal{T}$, which show that the expected value of $p_t(\rho,\rho)$ is controlled by $t^{-2/3}$ with at most $O(1)$ fluctuations as $t\rightarrow 0$.

{\propn\label{annealedheatkernelstate} Let $\rho$ be the root of $\mathcal{T}$, then there exist constants $C_{15}, C_{16}$ such that
$$C_{15} t^{-2/3}\leq\mathbf{E}\left(p_t(\rho,\rho)\right)\leq C_{16} t^{-2/3},\hspace{20pt}\forall t\in(0,1).$$}

At this point, a comparison with the results obtained by Barlow and Kumagai for the random walk on the incipient cluster for critical percolation on a regular tree, \cite{BarKum}, is pertinent. First, observe that the incipient infinite cluster can be constructed as a particular branching process conditioned to never become extinct and the self-similar CRT (see \cite{Aldous2}) can be constructed as the scaling limit of a similar branching process. Note also that the objects studied here and by Barlow and Kumagai are both measure-metric space trees and so similar probabilistic and analytic techniques for estimating the heat kernel may be applied to them. Consequently, it is not surprising that the quenched local heat kernel bounds of \cite{BarKum} exhibit log-logarithmic differences similar to those obtained in this article and furthermore, the annealed heat kernel behaviour at the root is also shown to be the same in both settings. It should be noted though that the volume bounds which are crucial for obtaining these heat kernel bounds are proved in very different ways. Here we use Brownian excursion properties, whereas in \cite{BarKum}, branching process arguments are applied. Unlike in \cite{BarKum}, we do not prove annealed off-diagonal heat kernel bounds. This is primarily because there is no canonical way of labeling vertices (apart from the root) in the CRT.

A more concrete connection between the Brownian motion on the CRT and simple random walks on random graph trees is provided in \cite{Croydoncbp}, where it is shown that if $\mathcal{T}_n$ is a family of random graph trees that, when rescaled, converge to the CRT in a suitable space, then the associated random walks $X^n$, when rescaled, converge to the Brownian motion on $\mathcal{T}$. In turn, this result is closely linked to the convergence of the rescaled height process (which measures the distance of a simple random walk from the root) on the incipient cluster for critical percolation on a regular tree, see \cite{Kesten}, providing further evidence in support of the fact that the observations of the previous paragraph are not merely coincidental.

The article is organised as follows. In Section \ref{prelim}, we provide definitions of and give a brief introduction to the main objects of the discussion, namely the normalised Brownian excursion and the CRT. In Section \ref{excursionprop} we prove several results about the sample paths of the normalised Brownian excursion, which are then used to obtain quenched and annealed volume bounds in Section \ref{upperbound}. In Section \ref{hkbounds}, we explain how already established results about dendrites and measure-metric spaces may be applied to the CRT to construct a process on $\mathcal{T}$ and derive the quenched transition density estimates of Theorems \ref{heatkernel}, \ref{globalheatkernelstate} and \ref{localheatkernelstate}. The annealed heat kernel bounds are proved in Section \ref{annealedheatkernelsection}. The proofs of both the quenched and annealed heat kernel bounds use the results established in Section \ref{upcrossings} about the upcrossings of the normalised Brownian excursion. Finally, we show that the process we have constructed is actually Brownian motion on $\mathcal{T}$ in Section \ref{BM}.

\section{Preliminaries}\label{prelim}
\setcounter{thm}{0}

\subsection{Normalised Brownian excursion}\label{brownianexcursion}

An important part of the definition of the CRT is the Brownian excursion, normalised to have length 1, a process which we will denote $W$. In this section, we provide two characterisations of the law of the Brownian excursion and describe the appropriate normalisation.

We begin by defining the space of excursions, $U$. First, let $U'$ be the space of functions $f:\mathbb{R}_+\rightarrow \mathbb{R}_+$ for which there exists a $\tau(f)\in(0,\infty)$ such that
$$f(t)>0\hspace{20pt}\Leftrightarrow\hspace{20pt}t\in (0,\tau(f)).$$
We shall take $U:=U'\cap C(\mathbb{R}_+,\mathbb{R}_+)$, the restriction to the continuous functions contained in $U'$. The space of excursions of length $s$ is then defined to be $U^{(s)}:=\{f\in U:\:\tau(f)=s\}$.

Our first description of the law of $W$ involves conditioning the It\^{o} excursion law, which arises from the Poisson process of excursions of a standard Brownian motion. We will denote by $B=(B_t)_{t\geq 0}$ a standard, 1-dimensional Brownian motion starting from 0, built on a probability space with probability measure $\mathbf{P}$. Since the It\^o excursion law has been widely studied, we shall omit most of the technicalities here. For more details of excursion laws for Markov processes, the reader is referred to \cite{RevuzYor}, Chapter XII, and \cite{RogWil}, Chapter VI.

Let $L_t$ be the local time of $B$ at 0, and $L_t^{-1}:=\inf\{s>0:\:L_s>t\}$ be its right continuous inverse. Wherever $L_{t^-}^{-1}\neq L_t^{-1}$, we define $e_t\in U$ to be the (positive) excursion at local time $t$. In particular,
$$e_t(s):=\left\{\begin{array}{ll}\left|B_{L_{t-}^{-1}+s}\right|,&\hspace{20pt}0\leq s\leq L_{t}^{-1}-L_{t^-}^{-1},\\0,&\hspace{20pt}s>L_{t}^{-1}-L_{t^-}^{-1}.\end{array}\right.$$
The set of excursions of $B$ is denoted by $\Pi:=\{(t,e_t) :\: L^{-1}_{t^-}\neq L^{-1}_t\}$. The key idea is that $\Pi$ is a Poisson process on $(0,\infty)\times U$. More specifically, there exists a $\sigma$-finite measure, $n$, on $U$ such that, under $\mathbf{P}$,
$$\# \left(\Pi \cap \cdot\right) {\buildrel d \over =} N(\cdot),$$
where $N$ is a Poisson random measure on $(0,\infty)\times U$ with intensity $dt\:n(df)$. Bearing this result in mind, even though it has infinite mass, the measure $n$ can be considered to be the ``law'' of the (unconditional) Brownian excursion.

Our second description of the law of the Brownian excursion is the well-known Bismut decomposition. In short, this characterisation of $n$ describes how, if $(t,f)$ is chosen according to the measure $\mathbf{1}_{[0,\tau(f)]}(t)dt\:n(df)$, then the law of $f(t)$ is Lebesgue measure on $\mathbb{R}_+$ and, conditionally on $f(t)=a$, the processes $(f(s\wedge t))_{s\geq 0}$ and $(f((\tau(f)-s)\vee t))_{s\geq 0}$ are independent 3-dimensional Bessel processes run until they last hit $a$. More precisely, let $m$ be the law of the 3-dimensional Bessel process, and note that, for every $a\geq0$ and $m$-a.e. $f\in C(\mathbb{R}_+,\mathbb{R}_+)$, the last hitting time of $a$ by $f$, defined by $T_a(f):=\sup\{t\geq 0:\:f(t)=a\}$ is finite. The result of interest may now be written as follows: for every non-negative measurable function $F$ on $C(\mathbb{R}_+,\mathbb{R}_+)^2$,
\begin{eqnarray*}
\lefteqn{\int n(df) \int_0^{\tau(f)} dt F((f(s\wedge t))_{s\geq 0},(f((\tau(f)-s)\vee t))_{s\geq 0})}\\
&=&2\int_0^\infty da \int m(df_1)\int m(df_2)F((f_1(s\wedge T_a(f_1)))_{s\geq 0},(f_2(s\wedge T_a(f_2)))_{s\geq 0}).
\end{eqnarray*}
See \cite{RevuzYor}, Theorem XII.4.7 for a proof.

We now describe how the excursion measure $n$ can be decomposed into a collection of related measures on the spaces of excursions of a fixed length. For $c>0$, the re-normalisation operator $\Lambda_c:U\rightarrow U$ is defined by
$$\Lambda_c(f)(t)=\frac{1}{\sqrt{c}}f(ct),\hspace{20pt}\forall t\geq 0,\:f\in U.$$
Clearly, if $f\in U$, then $\Lambda_{\tau(f)}(f)\in U^{(1)}$. Now, according to It\^o's description of $n$, (see \cite{RevuzYor}, Theorem XII.4.2), we can write, for any measurable $A\subseteq U$,
\begin{equation}\label{itodecomp}
n(A)=\int_0^\infty n^{(s)}(A\cap U^{(s)})\frac{ds}{\sqrt{2\pi s^3}},
\end{equation}
where, for each $s$, $n^{(s)}$ is a probability measure on $U^{(s)}$ that satisfies $n^{(s)}=n^{(1)}\circ \Lambda_s$, for some unique probability measure $n^{(1)}$ on $U^{(1)}$. In particular, we have that $n^{(1)}=n(\cdot|\:\tau(f)=1)$. A $U^{(1)}$-valued process which has law $n^{(1)}$ is said to be a normalised Brownian excursion. Henceforth, we assume that $W$ is a normalised Brownian excursion built on the probability space with probability measure $\mathbf{P}$.

\subsection{Continuum random tree}\label{crt}

The connection between trees and excursions is an area that has been of much recent interest. In this section, we look to provide a brief introduction to this link and also a definition of the CRT, which is the object of interest of this article.

Given a function $f\in U$, we define a distance on $[0,\tau(f)]$ by setting
\begin{equation}\label{distance}
d_f(s,t):=f(s)+f(t)-2m_f(s,t),
\end{equation}
where $m_f(s,t):=\inf\{f(r):\:r\in[s\wedge t,s\vee t]\}$. Then, we use the equivalence
$$s\sim t\hspace{20pt}\Leftrightarrow\hspace{20pt}d_f(s,t)=0,$$
to define $\mathcal{T}_f:=[0,\tau(f)]/\sim$. We can write this as $\mathcal{T}_f=\{\sigma_s:\:s\in[0,\tau(f)]\}$, where $\sigma_s:=[s]$ is the equivalence class containing $s$. It is then elementary (see \cite{LegallDuquesne}, Section 2) to check that
$$d_{\mathcal{T}_f}(\sigma_s,\sigma_t):=d_f(s,t),$$
defines a metric on $\mathcal{T}_f$, and also that $(\mathcal{T}_f,d_{\mathcal{T}_f})$ is a compact real tree in the sense of \cite{LegallDuquesne}, Definition 2.1 (there the term used is $\mathbb{R}$-tree). In particular, $\mathcal{T}_f$ is a dendrite: an arc-wise connected topological space, containing no subset homeomorphic to the circle. Furthermore, the metric $d_{\mathcal{T}_f}$ is a shortest path metric on $\mathcal{T}_f$, which means that it is additive along (injective) paths of $\mathcal{T}_f$. The root of the tree $\mathcal{T}_f$ is defined to be the equivalence class $\sigma_0$ and is denoted by $\rho_f$.

A natural volume measure to put on $\mathcal{T}_f$ is the projection of Lebesgue measure on $[0,\tau(f)]$. For open $A\subseteq\mathcal{T}_f$, let
$$\mu_f(A):=\ell\left(\{t\in[0,\tau(f)]:\:\sigma_t\in A\}\right),$$
where, throughout this article, $\ell$ is the usual 1-dimensional Lebesgue measure. This defines a Borel measure on $(\mathcal{T}_f,d_{\mathcal{T}_f})$, with total mass equal to $\tau(f)$.

The CRT is then simply the random dendrite that we get when the function $f$ is chosen according to the law of the normalised Brownian excursion. This differs from the Aldous CRT, which is based on the random function $2W$. Since this extra factor only has the effect of increasing distances by a factor of 2, our results are easily adapted to apply to Aldous' tree. In keeping with the notation used so far in this section, the measure-metric space should be written $(\mathcal{T}_W, d_{\mathcal{T}_W}, \mu_W)$, the distance on $[0,\tau(W)]$, defined at (\ref{distance}), $d_W$, and the root, $\rho_W$. However, we shall omit the subscripts $W$ with the understanding that we are discussing the CRT in this case. We note that $\tau(W)=1$, $\mathbf{P}$-a.s., and so $[0,\tau(W)]=[0,1]$ and $\mu$ is a probability measure on $\mathcal{T}$, $\mathbf{P}$-a.s. Finally, it is clear from the compactness of $\mathcal{T}$ that the diameter of $\mathcal{T}$, $\mathrm{diam}\mathcal{T}$, is finite $\mathbf{P}$-a.s.

\subsection{Other notation}

The $\delta$-level oscillations of a function $y$ on the interval $[s,t]$ will be written
$$\mathrm{osc}(y,[s,t],\delta):=\sup_{\buildrel{\scriptstyle{r,r'\in[s,t]:}}\over{\scriptstyle{|r'-r|\leq\delta}}}|y(r)-y(r')|,$$
and we will denote by $c_.$ constants taking a value in $(0,\infty)$. We shall also continue to use the notation introduced in Theorems \ref{globalstate} and \ref{heatkernel}, $\ln_1 x:=\ln x\vee 1$.

\section{Brownian excursion properties}\label{excursionprop}
\setcounter{thm}{0}

In this section, we use sample path properties of standard 1-dimensional Brownian motion and the 3-dimensional Bessel process to deduce various sample path properties for the normalised Brownian excursion. The definitions of the random variables $B$ and $W$ and measures $n$ and $n^{(1)}$ should be recalled from Section \ref{brownianexcursion}.

{\lem \label{modcont} $\mathbf{P}$-a.s.,
$$ \limsup_{\delta\rightarrow 0} \frac{\mathrm{osc} (W, [0,1], \delta)}{\sqrt{\delta\ln\delta^{-1}}}=\sqrt{2}.$$}
\begin{proof}
The proof of this lemma from L\'evy's 1937 result on the modulus of continuity of a standard Brownian motion in $\mathbb{R}$ (see \cite{RevuzYor}, Theorem I.2.7, for example) and the Poisson process description of the excursion law $n$ is standard, and so we omit it.
\end{proof}

{\lem \label{inflem} $\mathbf{P}$-a.s.,
\begin{equation}\label{rather}
\limsup_{\delta\rightarrow 0} \frac{\inf_{t\in[0,1-\delta]}\mathrm{osc} (W, [t,t+\delta], \delta)}{\sqrt{\delta(\ln\delta^{-1})^{-1}}}<\infty.
\end{equation}}
\begin{proof}
We start by proving the corresponding result for a standard Brownian motion. Fix a constant $c_{1}$ and then, for $n\geq 0$,
\begin{eqnarray*}
\lefteqn{\mathbf{P}\left(\inf_{j=0,\dots,2^n-1}\mathrm{osc}(B,[j2^{-n},(j+1)2^{-n}],2^{-n})\geq c_{1}2^{-\frac{n}{2}}n^{-\frac{1}{2}}\right)}\\
&= & \mathbf{P}\left(\mathrm{osc}(B,[j2^{-n},(j+1)2^{-n}],2^{-n})\geq c_{1}2^{-\frac{n}{2}}n^{-\frac{1}{2}},\:j=0,\dots,2^n-1\right)\\
&= & \mathbf{P}\left(\mathrm{osc}(B,[0,2^{-n}],2^{-n})\geq c_{1}2^{-\frac{n}{2}}n^{-\frac{1}{2}}\right)^{2^n},
\end{eqnarray*}
using the independent increments of a Brownian motion for the second equality. The probability in this expression is bounded above by
$$\mathbf{P}\left(\sup_{t\in[0,2^{-n}]}|B_t|\geq \frac{c_{1}}{2^{1+\frac{n}{2}}n^{\frac{1}{2}}}\right)= \mathbf{P}\left(T_{B_E(0,1)}\leq \frac{4n}{c_{1}^2}\right)\leq 1 - c_{2}e^{-\frac{c_{3}n}{c_{1}^2}},$$
for some constants $c_{2}, c_{3}$, where $T_{B_E(0,1)}$ represents the exit time of a standard Brownian motion from a Euclidean ball of radius 1 about the origin. The distribution of this random variable is known explicitly, see \cite{Cies}, and the above tail estimate is readily deduced from the expression that is given there. Using the fact that $1-x\leq e^{-x}$ for $x\geq0$ and summing over $n$, we have
$$\sum_{n=0}^{\infty}\mathbf{P}\left(\inf_{j=0,\dots,2^n-1}\mathrm{osc}(B,[j2^{-n},(j+1)2^{-n}],2^{-n})\geq c_{1}2^{-\frac{n}{2}}n^{-\frac{1}{2}}\right)\leq\sum_{n=0}^\infty e^{-2^n c_{2}e^{-\frac{c_{3}n}{c_{1}^2}}},$$
which is finite for $c_{1}$ chosen suitably large. Hence Borel-Cantelli implies that, $\mathbf{P}$-a.s., there exists a constant $c_{4}$ such that
$$\inf_{j=0,\dots,2^n-1}\mathrm{osc}(B,[j2^{-n},(j+1)2^{-n}],2^{-n})\leq c_{4} 2^{-\frac{n}{2}}n^{-\frac{1}{2}},\hspace{20pt}\forall n\geq 0.$$
Let $\delta\in(0,1]$, then $\delta\in[2^{-(n+1)},2^{-n}]$ for some $n\geq0$. Hence
\begin{eqnarray*}
\inf_{t\in[0,1-\delta]}\mathrm{osc}(B,[t,t+\delta],\delta)&\leq& \inf_{t\in[0,1-2^{-n}]}\mathrm{osc}(B,[t,t+2^{-n}],2^{-n})\\
&\leq& \inf_{j=0,\dots,2^n-1}\mathrm{osc}(B,[j2^{-n},(j+1)2^{-n}],2^{-n})\\
&\leq &c_{4} 2^{-\frac{n}{2}}n^{-\frac{1}{2}}\\
&\leq& c_{5} \sqrt{\delta(\ln\delta^{-1})^{-1}},
\end{eqnarray*}
which proves (\ref{rather}) holds when $W$ is replaced by $B$. By rescaling, an analogous result holds for any interval with rational endpoints. By countability and a monotonicity argument, this is easily extended to $\mathbf{P}$-a.s.,
$$\limsup_{\delta\rightarrow 0} \frac{\inf_{t\in[r,s-\delta]}\mathrm{osc} (B, [t,t+\delta], \delta)}{\sqrt{\delta(\ln\delta^{-1})^{-1}}}<\infty,\hspace{20pt}\forall 0\leq r<s<\infty.$$
The lemma can be deduced from this by considering the Poisson process of excursions of the Brownian motion $B$, as described in Section \ref{brownianexcursion}, and rescaling.
\end{proof}

{\lem \label{closeto0int} $\mathbf{P}$-a.s.,
$$ \limsup_{\delta\rightarrow 0}\frac{\int_0^1 \mathbf{1}_{\{W_t<\delta\}}dt}{\delta^2\ln\ln\delta^{-1}}=\frac{8}{\pi^2},$$
and
$$ \liminf_{\delta\rightarrow 0}\frac{\int_0^1 \mathbf{1}_{\{W_t<\delta\}}dt}{\delta^2\ln\ln\delta^{-1}}=2.$$}
\begin{proof} By applying the Bismut decomposition of $n$, see Section \ref{brownianexcursion}, to complete the proof, it will be sufficient to check that the lemma holds when $\int_0^1 \mathbf{1}_{\{W_t<\delta\}}dt$ is replaced by
\begin{equation}\label{occtime}
\int_0^\infty \mathbf{1}_{\{R^1_t<\delta\}}dt+\int_0^\infty \mathbf{1}_{\{R^2_t<\delta\}}dt,
\end{equation}
where, under the probability measure $\mathbf{P}$, the processes $R^1$ and $R^2$ are independent 3-dimensional Bessel processes. In the remainder of the proof we will write the first and second summands of (\ref{occtime}) as $T^1(\delta)$ and $T^2(\delta)$ respectively, and also $T(\delta):=T^1(\delta)+T^2(\delta)$.

The $\limsup$ asymptotic behaviour of the occupation time of 3-dimensional Brownian motion in a ball was considered in \cite{Cies}, Theorem 3, where it was shown that, $\mathbf{P}$-a.s.,
$$\limsup_{\delta\rightarrow 0}\frac{T^1(\delta)}{\delta^2\ln\ln\delta^{-1}}=\frac{8}{\pi^2}.$$
From this result we clearly have $\limsup_{\delta\rightarrow 0}T(\delta)/\delta^2\ln\ln\delta^{-1}\geq 8/\pi^2$. The identical upper bound will follow from a simple Borel-Cantelli argument, in a similar manner to the proof of \cite{Cies}, Theorem 3, if we can show that for every $\varepsilon>0$, there exists a constant $c_6$ such that
\begin{equation}\label{jod}
\mathbf{P}(T(\delta)>\delta^2 t)\leq c_6 e^{-(\frac{\pi^2}{8}-\varepsilon)t},
\end{equation}
for every strictly positive $t$ and $\delta$. To prove this first note that from \cite{Cies}, Theorem 1, we have $\mathbf{P}(T^1(\delta)>\delta^2 t)\leq e^{-\pi^2t/8}$ for every $t,\delta>0$. Consequently, using an argument analogous to the proof of \cite{TT}, Lemma 6.5, we are able to obtain, for $k\in\mathbb{N}$,
\begin{eqnarray*}
\mathbf{P}(T(\delta)> \delta^2 t)&=&\int_0^{\delta^2 t} \mathbf{P}(T^1(\delta)>\delta^2 t-s)d\mathbf{P}(T^1(\delta)\leq s)+\mathbf{P}(T^1(\delta)>\delta^2 t)\\
&=&\sum_{i=1}^k\int_{(i-1)\delta^2t/k}^{i\delta^2 t/k} \mathbf{P}(T^1(\delta)>\delta^2 t-s)d\mathbf{P}(T^1(\delta)\leq s)+\mathbf{P}(T^1(\delta)>\delta^2 t)\\
&\leq & \sum_{i=1}^k\mathbf{P}(T^1(\delta)>\delta^2 t(1-i/k))\mathbf{P}(T^1(\delta)>\delta^2 t(i-1)/k)\\
&&\hspace{220pt}+\:\mathbf{P}(T^1(\delta)>\delta^2 t)\\
&\leq& (k+1)e^{-\pi^2t(1-1/k)/8}.
\end{eqnarray*}
By choosing $k$ suitably large, we have the upper bound at (\ref{jod}), and the $\limsup$ result follows. For the $\liminf$ result, we can refer directly to \cite{TT}, Theorem 6.8, where it was shown that $\liminf_{\delta\rightarrow 0}T(\delta)/\delta^2(\ln\ln\delta^{-1})^{-1}=2$, $\mathbf{P}$-a.s.
\end{proof}

\section{Volume results}\label{upperbound}
\setcounter{thm}{0}

In this section, we apply properties of the normalised Brownian excursion, including those introduced in Section \ref{excursionprop}, to deduce the results about the volume growth on the CRT that were stated in the introduction. We start by proving the annealed volume result of Theorem \ref{annealedstate}; this follows easily from the expected occupation time of $[0,r)$ for a normalised Brownian excursion, for which an explicit expression is known.

\hspace{10pt}\\
{\bf Proof of Theorem \ref{annealedstate}:} By definition, we have that
\begin{equation}\label{rootexpress}
\mu(B(\rho,r))=\int_0^1 \mathbf{1}_{\{W_s<r\}}ds.
\end{equation}
An expression for the expectation of this random variable is obtained in \cite{DurIgl}, Section 3, giving
$$\mathbf{E}\left(\mu(B(\rho,r))\right)=\int_0^r 4ae^{-2a^2}da.$$
This integral is easily evaluated to give the desired result.
{\hspace*{\fill}$\square$\endlist\par\bigskip}

We now proceed to deduce the global upper volume bound for $\mathcal{T}$. The three main ingredients in the proof are the modulus of continuity result proved in Lemma \ref{modcont}, a bound on the tail of the distribution of the volume of a ball about the root, and the invariance under re-rooting of the CRT. This final property is also important in proving the local volume growth results. Before continuing, we define precisely what we mean by re-rooting and state the invariance result that we will use.

Given $W$ and $s\in[0,1]$, we define the shifted process $W^{(s)}=(W^{(s)}_t)_{0\leq t \leq 1}$ by
$$
W^{(s)}_t:= \left\{\begin{array}{ll}W_s+W_{s+t}-2m(s,s+t),&\hspace{20pt}0\leq t \leq 1-s\\ W_s+W_{s+t-1}-2m(s+t-1,s),&\hspace{20pt}1-s\leq t \leq 1,\end{array}\right.$$
where $m=m_W$ is the minimum function defined in Section \ref{crt}. It is known that, for any fixed $s$, the process $W^{(s)}$ has the same distribution as $W$, (see \cite{marmok}, Proposition 4.9). From this fact, an elementary application of Fubini's theorem allows it to be deduced that, for any measurable $A\subseteq U^{(1)}$,
\begin{equation}\label{shift}
\mathbf{E}(\ell\{s\in[0,1]:\:W^{(s)}\in A\})=\mathbf{P}(W\in A),
\end{equation}
where, as in Section \ref{crt}, $\ell$ is Lebesgue measure on [0,1]. Written down in terms of excursions, it is not immediately clear what this result is telling us about the CRT. However, if we observe that the real tree associated with $W^{(s)}$ has exactly the same structure as $\mathcal{T}$, apart from the fact that it is rooted at $\sigma_s$ instead of $\rho$, and recall that $\mu$ is the image of $\ell$ under the map $s\mapsto\sigma_s$, then we are able to obtain to conclude from (\ref{shift}) that, heuristically, any property of the random real tree $\mathcal{T}$ that holds when the root is assumed to be $\rho$ will also hold when we consider the root to be $\sigma$, for $\mu$-a.e. $\sigma\in\mathcal{T}$. More precisely, if we use the notation $\mathcal{T}^\sigma$ to represent the real tree $\mathcal{T}$ with root $\sigma\in\mathcal{T}$, then
\begin{equation}\label{reroot}
\mathbf{E}(\mu\{\sigma\in\mathcal{T}:\:\mathcal{T}^\sigma\in A\})=\mathbf{P}(\mathcal{T}\in A),
\end{equation}
for any measurable subset, $A$, of the space of compact rooted real trees.

The following result provides an exponential tail bound for the distribution of the volume of a ball of radius $r$ about the root that will be useful in the proof of Proposition \ref{globalsupupper}.

{\lem \label{tailbounds0} There exist constants $c_{7},c_{8}$ such that, for all $r>0,\lambda\geq 1$,
$$\mathbf{P}\left(\mu(B(\rho, r))\geq r^2\lambda\right)\leq c_{7}e^{-c_{8}\lambda}.$$}
\begin{proof} First, observe that we can consider $n^{(1)}$ to the law of the excursion of Brownian motion straddling a fixed time, conditional on this excursion having length 1, see \cite{RevuzYor}, Proposition XII.3.4. Hence it follows from \cite{Chung}, Theorem 9 that, for $r>0$,
$$\mathbf{E}\left(\left(\int_0^1\mathbf{1}_{\{W_t<r\}}dt\right)^k\right)\leq (k+1)!r^{2k}.$$
Recalling the expression for $\mu(B(\rho,r))$ that was stated at (\ref{rootexpress}) and applying an exponential version of Markov's inequality, we obtain from this that, for all $\theta\in(0,1)$,
$$\mathbf{P}\left(\mu(B(\rho,r))\geq r^2\lambda \right)\leq \mathbf{E}\left(e^{\theta r^{-2}\int_0^1\mathbf{1}_{\{W_t<r\}}dt-\theta\lambda}\right)\leq\sum_{k=0}^\infty\frac{\theta^k(k+1)!e^{-\theta\lambda}}{k!}=\frac{e^{-\theta\lambda}}{(1-\theta)^2}.$$
\end{proof}

{\propn \label{globalsupupper} $\mathbf{P}$-a.s., there exists a constant $c_{9}$ such that
$$\sup_{\sigma \in \mathcal{T}} \mu (B(\sigma, r))\leq c_{9} r^2 \ln_1 r^{-1},\hspace{20pt}\forall r\in (0,{\rm diam}\mathcal{T}).$$}
\begin{proof} In the proof, we shall denote $r_n:=e^{-n}$, $\delta_n:=r_n^2(\ln_1 r_n^{-1})^{-1}$, and also write $g(r):=r^2\ln_1r^{-1}$. Furthermore, we introduce the notation, for $\lambda\in(0,1]$, $n_0\in\mathbb{N}$,
$$A_\lambda(n_0):=\left\{\mathrm{osc}(W,[0,1],\lambda\delta_n)\leq r_n,\:\forall n\geq n_0\right\},$$
which represents a collection of sample paths of $W$ which are suitably regular for our purposes. We will start by showing that the claim holds on each set of the form $A_\lambda(n_0)$. Until otherwise stated, we shall assume that $\lambda$ and $n_0$ are fixed. Now, consider the sets
$$B_n:=\left\{\sup_{\sigma\in\mathcal{T}}\mu(B(\sigma,r_n))> c_{10} g(r_n)\right\}\cap A_\lambda(n_0),$$
where $n\geq n_0$, and $c_{10}$ is a constant we will specify below. Clearly, on the event $B_n$, the random subset of $[0,1]$ defined by
$$\mathcal{I}_n:=\left\{t\in[0,1]:\:|t-s|<\lambda\delta_n\mbox{ for some $s\in[0,1]$ with $\mu(B(\sigma_s,r_n))\geq c_{10} g(r_n)$}\right\}$$
has Lebesgue measure no less than $\lambda\delta_n$. Thus, if $U$ is a $U[0,1]$ random variable, independent of $W$, then
$$\mathbf{P}(U\in\mathcal{I}_n,\:B_n)=\mathbf{E}(\mathbf{P}(U\in\mathcal{I}_n|W){1}_{B_n})\geq\lambda\delta_n\mathbf{P}(B_n).$$
Moreover, on the event $\{U\in\mathcal{I}_n\}\cap B_n$, there exists an $s\in[0,1]$ for which both $|U-s|<\lambda\delta_n$ and $\mu(B(\sigma_s,r_n))\geq c_{10} g(r_n)$ are satisfied. Applying the modulus of continuity property that holds on $A_\lambda(n_0)$, for this $s$ we have that $d_{\mathcal{T}}(\sigma_s,\sigma_U)\leq 3r_n$, and so $\mu(B(\sigma_U,4r_n))\geq c_{10} g(r_n)$. Hence the above inequality implies that
\begin{eqnarray*}
\lambda\delta_n\mathbf{P}(B_n)&\leq&\mathbf{P}\left(\mu(B(\sigma_U,4r_n))\geq c_{10} g(r_n)\right)\\
&=&\mathbf{E}(\mu\{\sigma\in\mathcal{T}:\mu(B(\sigma,4r_n))\geq c_{10} g(r_n)\})\\
&=&\mathbf{P}\left(\mu(B(\rho,4r_n))\geq c_{10} g(r_n)\right)\\
&\leq&c_7e^{-c_8c_{10}n/16},
\end{eqnarray*}
Here, the first equality is a simple consequence of the fact $\mu$ is the image of Lebesgue measure on [0,1] under the map $s\mapsto\sigma_s$, we have applied the re-rooting result of (\ref{reroot}) to deduce the second equality, and the distributional tail bound of Lemma \ref{tailbounds0} is used to obtain the final line. As a consequence, we have that
$$\sum_{n\geq n_0}\mathbf{P}(B_n)\leq c_7\lambda^{-1}\sum_{n\geq n_0}ne^{2n}e^{-c_8c_{10}n/16},$$
which is finite for $c_{10}$ chosen suitably large. Appealing to Borel-Cantelli and applying a simple monotonicity argument, this implies that $\mathbf{P}$-a.s. on $A_\lambda(n_0)$
$$\limsup_{r \rightarrow 0}\frac{\sup_{\sigma\in\mathcal{T}}\mu(B(\sigma,r))}{r^2\ln_1r^{-1}}<\infty.$$
By countability, the same must be true on the set $A_\lambda:=\cup_{n\geq 1}A_\lambda(n)$. To deduce the proposition, note that Lemma \ref{modcont} gives $\mathbf{P}(A_{\lambda})=1$ if $\lambda$ is small enough.
\end{proof}

We continue by proving the global lower bounds for the volume of balls of the CRT. The estimates are simple corollaries of the excursion modulus of continuity results proved in Lemmas \ref{modcont} and \ref{inflem}.

{\propn \label{globallowerprop} $\mathbf{P}$-a.s., there exists a constant $c_{11}$ such that
$$\inf_{\sigma \in \mathcal{T}} \mu (B(\sigma, r))\geq c_{11} r^2 \left(\ln_1 r^{-1}\right)^{-1},\hspace{20pt}\forall r\in (0,{\rm diam}\mathcal{T}).$$}
\begin{proof}
For $s\in[0,1], r>0$, define
$$\alpha_u(s,r):=\inf\{t\geq 0:\:|W_{s+t}-W_s|> r\},$$
$$\alpha_l(s,r):=\inf\{t\geq 0:\:|W_s-W_{s-t}|> r\},$$
where these expressions are defined to be $1-s$, $s$ if the infimum is taken over an empty set, respectively. From the definition of $d$ at (\ref{distance}), it is readily deduced that $d(s,t)\leq r$ for all $t\in (s-\alpha_l(s,r/3),s+\alpha_u(s,r/3))$, from which it follows that $d_{\mathcal{T}}(\sigma_s,\sigma_t)\leq r$ for all t in this range. Hence
\begin{equation}\label{measurealpha}
\mu(B(\sigma_s,r))\geq \alpha_l(s,r/3)+\alpha_u(s,r/3).
\end{equation}

We now show how the right hand side of this inequality can be bounded below, uniformly in $s$, using the uniform modulus of continuity of the Brownian excursion. By Lemma \ref{modcont}, $\mathbf{P}$-a.s. there exist constants $c_{12}, \eta\in(0,1)$ such that
\begin{equation}\label{oscbound}
\mathrm{osc}(W,[0,1],\delta)\leq c_{12} \sqrt{\delta \ln\delta^{-1}},\hspace{20pt}\forall \delta\in(0,\eta).
\end{equation}
Set $r_{1}=3 c_{12} \sqrt{\eta \ln\eta^{-1}}$ and then, for $r\in(0,r_{1})$, choose $\delta=\delta(r)$ to satisfy $r=3 c_{12} \sqrt{\delta \ln\delta^{-1}}$. It follows from the inequality at (\ref{oscbound}) that if $r\in(0,r_{1})$ and $|W_s-W_t|>r/3$, then $|s-t|>\delta(r)\geq c_{13} r^2 \left(\ln r^{-1}\right)^{-1}$,
where $c_{13}$ is a constant depending only on $c_{12}$ and $\eta$. By definition, this implies that
$$\alpha_l(s,r/3)\geq c_{13} r^2 \left(\ln r^{-1}\right)^{-1}\wedge s,\hspace{20pt}\alpha_u(s,r/3)\geq c_{13} r^2 \left(\ln r^{-1}\right)^{-1}\wedge (1-s),$$
for all $s\in[0,1]$, $r\in(0,r_{1})$. Adding these two expressions and taking a suitably small constant, we obtain that $\mathbf{P}$-a.s., there exists a constant $c_{14}$ such that
$$\inf_{s\in[0,1]}\left(\alpha_l(s,r/3)+\alpha_u(s,r/3)\right)\geq c_{14} r^2 \left(\ln_1 r^{-1}\right)^{-1},\hspace{20pt}\forall r \in (0,\mathrm{diam}\mathcal{T}),$$
where we use the fact that $\mathrm{diam}\mathcal{T}$ is $\mathbf{P}$-a.s. finite. Taking infima in (\ref{measurealpha}) and comparing with the above inequality completes the proof.
\end{proof}

{\propn \label{globalsuplower} $\mathbf{P}$-a.s., there exists a constant $c_{15}$ such that
$$\sup_{\sigma \in \mathcal{T}} \mu (B(\sigma, r))\geq c_{15} r^2 \ln_1 r^{-1},\hspace{20pt}\forall r\in (0,{\rm diam}\mathcal{T}).$$}
\begin{proof}
By following an argument similar to that used in the proof of the previous proposition to transfer the excursion result to a result about the volume of balls in the CRT, the proposition may be deduced from Lemma \ref{inflem}.
\end{proof}

We conclude this section by discussing the proofs of the remaining volume growth results of the introduction. Firstly, given the characterisation of the volume of a ball of radius $r$ around the root that appears at (\ref{rootexpress}) and the re-rooting result of (\ref{reroot}), the local volume growth asymptotics of Theorem \ref{localstate} are a straightforward consequence of the result proved in Lemma \ref{closeto0int} about the time spent close to 0 by a normalised Brownian excursion. Finally, for the proof of Theorem \ref{globalstate}, note that Propositions \ref{globalsupupper} and \ref{globalsuplower} contain the upper and lower bounds for $\sup\mu(B(\sigma,r))$ respectively. The lower bound for $\inf\mu(B(\sigma,r))$ was proved in Proposition \ref{globallowerprop}; the upper bound follows easily from the local results.

\section{Brownian excursion upcrossings}\label{upcrossings}
\setcounter{thm}{0}

To complete the proofs of the heat kernel bounds in Sections \ref{hkbounds} and \ref{annealedheatkernelsection}, as well as the volume estimates we have already obtained, we need some extra information about the local structure of the CRT. This will follow from the results about the upcrossings of a normalised Brownian excursion that we prove in this section.

Henceforth, we define, for $f\in U$,
$$N_a^b(f) :=\#\{\mbox{upcrossings of $[a,b]$ by $f$}\}$$
and $N_a^b:=N_a^b(W)$ to be the number of upcrossings of $[a,b]$ by the normalised Brownian excursion. Also, for $f\in U$, define $h(f):=\sup\{f(t):\:t\geq 0\}$, the height of the excursion function. It is well known, (see \cite{RevuzYor}, Chapter XII) that the tail of the ``distribution'' of $h(f)$ under $n$ is simply given by
\begin{equation}\label{heighttail}
n(h(f)\geq x)=\frac{1}{x},\hspace{20pt}\forall x>0.
\end{equation}

We now calculate the generating function of $N_\delta^{2\delta}(f)$ under the probability measure $n(\cdot|\:h(f)\geq \delta)$. Notice that the expression we obtain does not depend on $\delta$.

{\lem \label{upcrossgen} For $z<2$, $\delta> 0$,
$$n(z^{N_\delta^{2\delta}(f)}|\:h(f)\geq\delta)=\frac{1}{2-z}.$$
}
\begin{proof}
Suppose $X=(X_t)_{t\geq 0}$ is a process with law $n(\cdot|\:h(f)\geq\delta)$, and define the stopping time $T:=\inf\{t>0:\:X_t=\delta\}$ to be the first time $X$ hits the level $\delta$. By the strong Markov property of $n$ (see \cite{RevuzYor}, Theorem XII.4.1), the process $X':=(X_{T+t})_{t\geq 0}$ is a standard Brownian motion started from $\delta$ and killed on hitting 0. Moreover, we clearly have $N_{\delta}^{2\delta}(X)=N_{\delta}^{2\delta}(X')$. An elementary argument using the symmetry and strong Markov property of linear Brownian motion implies that $N_{\delta}^{2\delta}(X')$ is a geometric, parameter $\frac{1}{2}$, random variable. The result follows.
\end{proof}

In the proof of the following result about the tail of the distribution of $N_\delta^{2\delta}(f)$, we shall use a result of Le Gall and Duquesne that states that the set of excursions above a fixed level form a certain Poisson process. We outline briefly the result here, full details may be found in \cite{LegallDuquesne}, Section 3. Fix $a>0$ and denote by $(\alpha_j,\beta_j)$, $j\in\mathcal{J}$, the connected components of the open set $\{s\geq 0:\:f(s)>a\}$. For any $j\in\mathcal{J}$, denote by $f^j$ the corresponding excursion of $f$ defined by:
$$f^j(s):=f((\alpha_j+s)\wedge\beta_j)-a,\hspace{20pt}s\geq 0,$$
and let $\tilde{f}^a$ represent the evolution of $f$ below the level $a$. More precisely, $\tilde{f}^a(s):=f(\tilde{\tau}^a_s)$, where
$$\tilde{\tau}^a_s:=\inf\left\{t\geq 0:\:\int_0^t\mathbf{1}_{\{f(r)\leq a\}}dr>s\right\}.$$
Applying the Poisson mapping theorem to \cite{LegallDuquesne}, Corollary 3.2, we find that under the probability measure $n(\cdot|\:h(f)>a)$ and conditionally on $\tilde{f}^a$, the point measure
$$\sum_{j\in\mathcal{J}} \delta_{f^j}(df')$$
is distributed as a Poisson random measure on $U$ with intensity given by a multiple of $n(df')$. The relevant scaling factor is given by the local time of $f$ at the level $a$, which we shall denote by $l^a$. Note that this is a measurable function of $\tilde{f}^a$.

{\lem \label{upcrossingbound} For $z\in[1,2)$, $\varepsilon>0$, $\delta\in (0,\varepsilon/2)$,
$$n(N_\delta^{2\delta}(f)\geq\lambda,\:h(f)\geq\varepsilon)\leq \frac{2z^{1-\lambda}}{(2-z)^2\varepsilon},\hspace{20pt}\forall\lambda\geq 0.$$}
\begin{proof}
For brevity, during this lemma we shall write $h=h(f)$ and $h^j=h(f^j)$, where $f^j$, $j\in\mathcal{J}$ are the excursions above the level $\delta$. Note that it is elementary to show that the quantity $N_\delta^{2\delta}(f)$ also counts the number of excursions of $f$ above the level $\delta$ which hit the level $2\delta$. Thus
\begin{eqnarray*}
N_{\delta}^{2\delta}(f)&=&\#\{j\in\mathcal{J}:\:h^j\geq\delta\}\\
&=&\#\{j\in\mathcal{J}:\:h^j\in[\delta,\varepsilon-\delta)\}+\#\{j\in\mathcal{J}:\:h^j\geq\varepsilon-\delta\}.
\end{eqnarray*}
We shall denote these two summands $N_1$ and $N_2$, respectively. From the observation preceding this lemma that the excursions above the level $\delta$ form a Poisson process on $U$ and the fact that the sets $\{h\in[\delta,\varepsilon-\delta)\}$ and $\{h\geq\varepsilon-\delta\}$ are disjoint, we can conclude that the $N_1$ and $N_2$ are independent under the measure $n(\cdot|\:h\geq\delta)$ and conditional on $\tilde{f}^\delta$. Furthermore, we note that $h\geq\varepsilon$ if and only if $N_2\geq1$. Hence, for $z\in(0,1]$,
\begin{eqnarray}
\lefteqn{n(z^{N_{\delta}^{2\delta}(f)}\mathbf{1}_{\{h\geq \varepsilon\}}\:|\:h\geq\delta)}\nonumber\\
&=&n(n(z^{N_1+N_2}\mathbf{1}_{\{N_2\geq 1\}}\:|\:\tilde{f}^\delta, h\geq \delta)\:|\:h\geq\delta)\nonumber\\
&=&n(n(z^{ N_1}\:|\:\tilde{f}^\delta, h\geq \delta)n(z^{ N_2}\mathbf{1}_{\{N_2\geq 1\}}\:|\:\tilde{f}^\delta, h\geq \delta)\:|\:h\geq\delta).\label{earlier}
\end{eqnarray}
Since on the set $\{h\geq\delta\}$, and conditional on $\tilde{f}^\delta$, $N_1$ and $N_2$ are Poisson random variables with means $l^\delta n(h\in[\delta,\varepsilon-\delta))$ and $l^\delta n(h\geq\varepsilon-\delta)$ respectively, it is elementary to conclude that
$$n(z^{ N_1}\:|\:\tilde{f}^\delta, h\geq \delta)=e^{-l^\delta n(h\in[\delta,\varepsilon-\delta))(1-z)},$$
and also
$$n(z^{ N_2}\mathbf{1}_{\{N_2\geq 1\}}\:|\:\tilde{f}^\delta, h\geq \delta)=e^{-l^\delta n(h\geq\varepsilon-\delta)(1-z)}(1-e^{-l^\delta n(h\geq\varepsilon-\delta)z}).$$
Substituting these expressions back into (\ref{earlier}) and applying the formula for the excursion height distribution that was stated at (\ref{heighttail}), it follows that
\begin{eqnarray}
n(z^{N_{\delta}^{2\delta}(f)}\mathbf{1}_{\{h\geq \varepsilon\}}\:|\:h\geq\delta)&=&n(e^{-l^\delta (1-z)\delta^{-1}}\:|\:h\geq\delta)\nonumber\\
& &\hspace{40pt}-n(e^{-l^\delta[(1-z)\delta^{-1}+z(\varepsilon-\delta)^{-1}]}\:|\:h\geq \delta).\label{notsoearly}
\end{eqnarray}
As in the proof of \cite{LeGallDuquesne2}, Theorem 1.4.1, for example, it is possible to use the strong Markov property of the excursion measure $n$ to deduce that $l^\delta$ satisfies $n(1-e^{-\lambda l^\delta})=\lambda(1+\lambda\delta)^{-1}$, for $\lambda\geq 0$. It follows that, under $n(\cdot\:|\:h\geq\delta)$, $l^\delta$ is an exponential, mean $\delta$, random variable. Returning to (\ref{notsoearly}), this fact implies that
$$n(z^{N_{\delta}^{2\delta}(f)}\mathbf{1}_{\{h\geq \varepsilon\}}\:|\:h\geq\delta)=\frac{1}{2-z}-\frac{1}{2-z+\delta z(\varepsilon-\delta)^{-1}}.$$
By a simple analytic continuation argument, we may extend the range of $z$ for which this holds to [0,2). Finally, for $z\in[1,2)$, we have that, because $\delta<\varepsilon/2$,
\begin{eqnarray*}
n(N_\delta^{2\delta}(f)\geq\lambda,\:h(f)\geq\varepsilon)&\leq&n(z^{N_{\delta}^{2\delta}(f)}1_{h\geq \varepsilon})z^{-\lambda}\\
&=&n(z^{N_{\delta}^{2\delta}(f)}1_{h\geq \varepsilon}\:|\:h\geq\delta)n(h\geq\delta)z^{-\lambda}\\
&=&\frac{1}{\delta z^{\lambda}}\left(\frac{1}{2-z}-\frac{1}{2-z+\delta z(\varepsilon-\delta)^{-1}}\right)\\
&\leq&\frac{2z^{1-\lambda}}{(2-z)^2\varepsilon},
\end{eqnarray*}
which completes the proof.
\end{proof}

We now reach the main results of this section, which give us an upper bound on the upcrossings of $[\delta,4\delta]$ by the normalised Brownian excursion for small $\delta$ and also a distribution tail estimate.

{\propn \label{upcross} $\mathbf{P}$-a.s.,
$$\limsup_{\delta\rightarrow 0} \frac{N_\delta^{4\delta}}{\ln\ln\delta^{-1}}<\infty.$$}
\begin{proof}
Fix $z \in (1,2)$, $\varepsilon\in(0,1)$ and let $\lambda$ be a constant. Then, taking $m_0$ suitably large, we have by the previous lemma that
$$\sum_{m=m_0}^\infty n(N_{2^{-m}}^{2^{1-m}}(f)\geq \lambda \ln m,\:h(f)\geq\varepsilon)\leq\sum_{m=m_0}^\infty \frac{2m^{(1-\lambda) \ln z}}{(2-z)^2\varepsilon} <\infty,$$
for $\lambda$ chosen suitably large. Hence an application of Borel-Cantelli implies that, on $\{h(f)>\varepsilon\}$,
$$\limsup_{m\rightarrow\infty} \frac{N_{2^{-m}}^{2^{1-m}}(f)}{\ln m} <\infty,\hspace{20pt}n\mbox{-a.s.}$$
Now for $\delta\in[2^{-(m+1)}, 2^{-m}]$, we have $[2^{-m},2^{1-m}]\subseteq[\delta,4\delta]$, and so $N_\delta^{4\delta}(f)\leq N_{2^{-m}}^{2^{1-m}}(f)$. Thus, on $\{h(f)>\varepsilon\}$,
$$\limsup_{\delta\rightarrow 0} \frac{N_{\delta}^{4\delta}(f)}{\ln \ln\delta^{-1}} \leq \limsup_{m\rightarrow\infty} \frac{N_{2^{-m}}^{2^{1-m}}(f)}{\ln m} <\infty,\hspace{20pt}n\mbox{-a.s.}$$
By $\sigma$-additivity, this is easily extended to hold on $U$, $n$-a.s. A simple rescaling argument allows us to obtain the same result $n^{(1)}$-a.s. on $U^{(1)}$. Since $n^{(1)}$ is the law of $W$, we are done.
\end{proof}

{\propn \label{upcrossingtail} There exist constants $c_{16}, c_{17}$ such that
$$\mathbf{P}\left(N_\delta^{4\delta}\geq \lambda\right) \leq c_{16}e^{-c_{17}\lambda},\hspace{20pt}\forall \delta>0,\lambda\geq 1.$$}
\begin{proof} Using It\^o's description of $n$, (\ref{itodecomp}), it is possible to deduce the following alternative characterisation of $n^{(1)}$. For measurable $A\subseteq U^{(1)}$, $n^{(1)}(A)=n(\Lambda_\tau(f)\in A|\:\tau\in[1,2])$.
Hence
$$\mathbf{P}(N_\delta^{4\delta}\geq \lambda)=n^{(1)}(N_\delta^{4\delta}(f)\geq \lambda)= c_{18}n(N_{\delta\sqrt\tau}^{4\delta\sqrt\tau}(f)\geq \lambda,\:\tau\in[1,2]).$$
However, for $\tau\in[1,2]$ we have that $[\sqrt 2 \delta,2\sqrt 2 \delta]\subseteq [\delta\sqrt\tau,4\delta\sqrt\tau]$, and so $N_{\delta\sqrt\tau}^{4\delta\sqrt\tau}(f)\leq N_{\sqrt 2 \delta}^{2\sqrt 2 \delta}(f)$. Hence
\begin{equation}\label{somekindofbound}
\mathbf{P}(N_\delta^{4\delta}\geq \lambda)\leq c_{18}n(N_{\sqrt 2 \delta}^{2\sqrt 2\delta}(f)\geq \lambda,\:\tau\in[1,2]).
\end{equation}
We now consider the cases $\lambda\leq \delta^{-1}$ and $\lambda\geq \delta^{-1}$ separately. First, suppose $\lambda\geq\delta^{-1}$. Since $\lambda\geq 1$, if $N_{\sqrt 2 \delta}^{2\sqrt 2 \delta}(f)\geq \lambda$, then $h(f)\geq \sqrt 2\delta$. Thus
\begin{eqnarray*}
\mathbf{P}(N_\delta^{4\delta}\geq \lambda)&\leq & c_{18} n(N_{\sqrt 2\delta}^{2\sqrt 2\delta}(f)\geq \lambda,\:h\geq \sqrt 2 \delta)\\
&=& \frac{c_{18}}{\sqrt2\delta} n\left(N_{\sqrt 2\delta}^{2\sqrt 2\delta}(f)\geq \lambda\:\vline\:h\geq \sqrt 2 \delta\right)\\
&\leq& \frac{c_{18}z^{-\lambda}}{\sqrt2\delta} n\left(z^{N_{\sqrt 2\delta}^{2\sqrt 2\delta}(f)}\:\vline\:h\geq \sqrt 2 \delta\right)\\
&\leq&\frac{c_{19}\lambda z^{-\lambda}}{2-z},
\end{eqnarray*}
for $z\in(1,2)$, by Lemma \ref{upcrossgen}. By fixing $z\in(1,2)$, this is clearly bounded above by $c_{20} e^{-c_{21}\lambda}$ for suitable choice of $c_{20},c_{21}$.

The case $\lambda\leq \delta^{-1}$ requires a little more work. We first break the upper bound of (\ref{somekindofbound}) in two parts. For $\varepsilon>0$,
\begin{equation}\label{somekindofbound2}
\mathbf{P}(N_\delta^{4\delta}\geq \lambda)\leq c_{18}n(N_{\sqrt 2 \delta}^{2\sqrt 2\delta}(f)\geq \lambda,\:h\geq \varepsilon)+c_{18}n(\tau\in[1,2],\:h<\varepsilon).
\end{equation}
An upper bound for the first term of the form $c_{22}e^{-c_{23}\lambda}\varepsilon^{-1}$ is given by Lemma \ref{upcrossingbound} when $\delta\leq\varepsilon/2\sqrt2$. For the second term, we again appeal to It\^o's description of $n$, (\ref{itodecomp}), to obtain
$$n(\tau\in[1,2],\:h<\varepsilon)=\int_1^2 n^{(s)}(h<\varepsilon)\frac{ds}{\sqrt{2\pi s^3}},$$
where $n^{(s)}$ is a probability measure that satisfies $n^{(s)}(A)=n^{(1)}(\Lambda_s(A))$. However, this scaling property implies that $n^{(s)}(h<\varepsilon)$ is decreasing in $s$. Consequently, we have that $n(\tau\in[1,2],\:h<\varepsilon)\leq (2\pi)^{-1/2} n^{(1)}(h<\varepsilon)$. The right hand side of this expression represents the distribution of the maximum of a normalised Brownian excursion, which is known exactly, see \cite{Chung}, Theorem 7. In particular, we have that
$$n^{(1)}(h<\varepsilon)=c_{24}\varepsilon^{-3}\sum_{m=1}^\infty m^2 e^{-\frac{m^2\pi^2}{2\varepsilon^2}},$$
and some elementary analysis shows this is bounded above by $c_{25}\varepsilon^{-3}e^{-\frac{c_{26}}{\varepsilon}}$. By taking $\varepsilon=2\sqrt2 \lambda^{-1}$, we can use these bounds on the first and second terms of (\ref{somekindofbound2}) to obtain the desired result.
\end{proof}

\section{Quenched heat kernel bounds}\label{hkbounds}
\setcounter{thm}{0}

The heat kernel estimates deduced in this section are a straightforward application of two main ideas. Firstly, we apply results of \cite{Kigamidendrite}, by Kigami, to construct a natural Dirichlet form on $(\mathcal{T}, d_{\mathcal{T}}, \mu)$. Associated with such a form is a Laplacian, $\Delta_\mathcal{T}$, on $\mathcal{T}$. Secondly, we use results of \cite{Croydon}, which show that the volume bounds we have already obtained are sufficient to deduce the existence of a heat kernel for $\Delta_\mathcal{T}$ and bounds on it.

The following result is proved by Kigami as Theorem 5.4 of \cite{Kigamidendrite}. Most of the terms used are standard, and so we will not define them here. Definition 0.5 of \cite{Kigamidendrite} specifies the precise conditions that make a symmetric, non-negative quadratic form a finite resistance form. For more examples of this type of form, see \cite{Kigami}. We shall not explain here how to construct the finite resistance form associated with a shortest path metric on a dendrite. Full details are given in Section 3 of \cite{Kigamidendrite}.

{\lem Let $K$ be a dendrite and let $d$ be a shortest path metric on $K$. Suppose $(K,d)$ is locally compact and complete, then $(\mathcal{E},\mathcal{F}\cap L^2(K,\nu))$, where $(\mathcal{E},\mathcal{F})$ is the finite resistance form associated with $(K,d)$, is a local, regular Dirichlet form on $L^2(K,\nu)$ for any $\sigma$-finite Borel measure $\nu$ on $K$ that satisfies $\nu(A)>0$ for every non-empty open set $A\subseteq K$.}

\hspace{10pt}

Consider $f\in U^{(1)}$. As remarked in Section \ref{crt}, $\mathcal{T}_f$ is a compact dendrite and $d_{\mathcal{T}_f}$ a shortest path metric on $\mathcal{T}_f$. Moreover, using simple path properties of the Brownian excursion, it is easy to check that $\mu_f$ satisfies the measure conditions of the above lemma for $f\in\tilde{U}$, where $\tilde{U}\subseteq U^{(1)}$ is a set which satisfies $\mathbf{P}(W\in \tilde{U})=1$. Consequently, we can define the finite resistance form $(\mathcal{E}_\mathcal{T},\mathcal{F}_\mathcal{T})$ associated with $(\mathcal{T},d_\mathcal{T})$ and  $(\mathcal{E}_\mathcal{T},\mathcal{F}_\mathcal{T}\cap L^2(\mathcal{T},\mu))$ is a local, regular Dirichlet form on $L^2(\mathcal{T},\mu)$, $\mathbf{P}$-a.s.

In fact, it is also proved in \cite{Kigamidendrite} that the correspondence between shortest path metrics on dendrites and resistance forms is one-to-one in a certain sense. Specifically, if $W\in \tilde{U}$ and so $(\mathcal{E}_\mathcal{T},\mathcal{F}_\mathcal{T})$ exists, define the resistance function by
\begin{equation}\label{resistance}
R(A,B)^{-1}:=\inf\left\{\mathcal{E}_\mathcal{T}(u,u):\:u\in\mathcal{F}_\mathcal{T},u|_A=1,u|_B=0\right\},
\end{equation}
for disjoint subsets $A, B$ of $\mathcal{T}$. We can recover $d_\mathcal{T}$ by taking, for $\sigma,\sigma'\in\mathcal{T}$, $\sigma\neq\sigma'$, $d_\mathcal{T}(\sigma,\sigma')=R(\{\sigma\},\{\sigma'\})$ and $d_\mathcal{T}(\sigma,\sigma)=0$. This means that the metric $d_\mathcal{T}$ is the effective resistance metric associated with $(\mathcal{E}_\mathcal{T},\mathcal{F}_\mathcal{T})$, see Corollary 3.4 of \cite{Kigamidendrite} for a proof of this. Using this description of $d_{\mathcal{T}}$ and a standard rescaling argument, it is easily deduced that\begin{equation}\label{contin}
(u(\sigma)-u(\sigma'))^2\leq \mathcal{E}(u,u)d_{\mathcal{T}}(\sigma,\sigma'),\hspace{20pt}\forall \sigma, \sigma'\in\mathcal{T},\:u\in\mathcal{F}_\mathcal{T},
\end{equation}
which has as a consequence that $\mathcal{F}_\mathcal{T}$ is contained in $C(\mathcal{T})$, the space of continuous functions on $(\mathcal{T},d_\mathcal{T})$. Since $\mathcal{T}$ is compact, we must also have $C(\mathcal{T})\subseteq L^2(\mathcal{T},\mu)$, meaning that the Dirichlet form of interest is simply $(\mathcal{E}_\mathcal{T},\mathcal{F}_\mathcal{T})$.

Given the Dirichlet form $(\mathcal{E}_\mathcal{T},\mathcal{F}_\mathcal{T})$, we can use the standard association to define a non-negative self-adjoint operator, $-\Delta_\mathcal{T}$, which has domain dense in $L^2(\mathcal{T},\mu)$ and satisfies
$$\mathcal{E}_\mathcal{T}(u,v)=-\int_\mathcal{T}u\Delta_{\mathcal{T}} vd\mu,\hspace{20pt}\forall u\in\mathcal{F}_\mathcal{T},v\in\mathcal{D}(\Delta_\mathcal{T}).$$
We can use this to define a reversible strong Markov process,
$$X=((X_t)_{t\geq 0}, \mathbf{P}_\sigma^{\mathcal{T}}, \sigma\in\mathcal{T}),$$with semi-group given by $P_t:=e^{t\Delta_\mathcal{T}}$. In fact, the locality of our Dirichlet form ensures that the process $X$ is a diffusion on $\mathcal{T}$.

As remarked in the introduction, a key factor in the description of the transition density of $X$, if it exists, is the volume growth of the space. The volume bounds we have already obtained for $\mathcal{T}$ mean that we can directly apply the bounds obtained in \cite{Croydon} for a resistance form with non-uniform volume doubling. Since stating the general theorem for resistance forms (see \cite{Croydon}, Theorem 4.1) would require a lot of extra notation, we omit to do so here and simply state the relevant conclusion for the CRT.

{\thm \label{process} $\mathbf{P}$-a.s., there is a reversible strong Markov diffusion $X$ on $\mathcal{T}$ with invariant measure $\mu$ and transition density $(p_t(\sigma,\sigma'))_{\sigma,\sigma'\in\mathcal{T},t>0}$ that satisfies the bounds at (\ref{lowerbound1}) and (\ref{upperbound1}).
}

\hspace{10pt}

Since, if it exists, the transition density of the process $X$ is a heat kernel of $\Delta_\mathcal{T}$, we can state the previous result in the following alternative form. Note, the main difference between a heat kernel and a transition density is that a heat kernel is defined only $\mu$-a.e., whereas the transition density is defined everywhere. Hence, for an arbitrary heat kernel of $\Delta_\mathcal{T}$, the bounds we have proved will hold only $\mu$-a.e. For full definitions of these two objects, see \cite{Croydon}, Section 5.

{\cor $\mathbf{P}$-a.s., there exists a local, regular Dirichlet form associated with the measure-metric space $(\mathcal{T},d_\mathcal{T},\mu)$. The related non-positive self-adjoint operator, $\Delta_\mathcal{T}$, admits a heat kernel $(p_t(\sigma,\sigma'))_{\sigma,\sigma'\in\mathcal{T},t>0}$ that satisfies the bounds at (\ref{lowerbound1}) and (\ref{upperbound1}).}

\hspace{10pt}

The proof of the remaining quenched heat kernel bounds also employ the techniques used in \cite{Croydon}. However, as well as the volume bounds, we need to apply the following extra fact about the asymptotics of the resistance from the centre of a ball to its surface, as $r\rightarrow 0$.

{\lem \label{resistance1} $\mathbf{P}$-a.s., for $\mu$-a.e. $\sigma\in\mathcal{T}$, there exist constants $c_{27}, r_{2}$ such that
$$c_{27}r\left(\ln\ln r^{-1}\right)^{-1}\leq R(\{\sigma\},B(\sigma,r)^c)\leq r,\hspace{20pt}\forall r\in(0,r_{2}).$$}
\begin{proof}
Choosing $r_{2}$ to be small enough so that $\sigma$ is connected to $B(\sigma,r)^c$ by a path of length $r$ immediately implies the upper bound. For the lower bound, following the argument of \cite{BarKum}, Lemma 4.4 we obtain
$$R(\{\sigma\},B(\sigma,r)^c)^{-1}\leq\frac{8M(\sigma,r)}{r},$$
where $M(\sigma,r)$ is defined to be the smallest number such that there exists a set $A=\{\sigma_1,\dots,\sigma_{M(\sigma,r)}\}$ with $d_{\mathcal{T}}(\sigma,\sigma_i)=r/4$ for each $i$, such that any path from $\sigma$ to $B(\sigma,r)^c$ must pass through the set $A$. For the CRT it is elementary to deduce that $M(\rho,r)\leq N_{r/4}^r$, where $N_{r/4}^r$ is the number of upcrossings of $[r/4,r]$ by $W$, as defined in Section \ref{upcrossings}. Thus, applying Proposition \ref{upcross}, the result holds at the root. This may be extended to hold $\mu$-a.e. using the re-rooting result of (\ref{reroot}).
\end{proof}

\hspace{10pt}\\
{\bf Proof of Theorem \ref{localheatkernelstate}:} On measure-metric spaces equipped with a resistance form, an upper bound for the on-diagonal part of the heat kernel of the form
\begin{equation}\label{heatkernelupper}
p_{2r\mu(B(\sigma, r))}(\sigma,\sigma)\leq\frac{2}{\mu(B(\sigma,r))}
\end{equation}
follows from a relatively simple analytic argument from the lower bound on the volume growth, see \cite{Kumagai}, Proposition 4.1 for an example of this kind of proof. Applying this, the two upper bounds of this result follow easily from the lower local volume bounds of Theorem \ref{localstate} and so we omit their proof here.

It remains to prove the lower local heat kernel bound. Again, the proof of this result is standard and so we shall only outline it briefly here. First, the resistance result of Lemma \ref{resistance1} and the local volume results of Theorem \ref{localstate} allow us to deduce that $\mathbf{P}$-a.s. for $\mu$-a.e. $\sigma\in\mathcal{T}$, there exist constants $c_{28}, c_{29}, r_{3}>0$ such that, for $r\in(0,r_{3})$,
$$\mathbf{E}_{\sigma'}^\mathcal{T} T_{B(\sigma,r)} \leq c_{28} r^3 \ln\ln r^{-1},\hspace{20pt}\forall \sigma'\in B(\sigma,r),$$
$$\mathbf{E}_{\sigma}^\mathcal{T} T_{B(\sigma,r)}\geq c_{29} r^3\left(\ln\ln r^{-1}\right)^{-4},$$
by applying an argument similar to the proof of \cite{Croydon}, Lemma 6.5. Here, $T_{B(\sigma,r)}$ is the exit time of the process $X$ from the ball $B(\sigma,r)$. Since $X$ is a Markov process we have that
$$\mathbf{E}_{\sigma}^\mathcal{T} T_{B(\sigma,r)}\leq t+\mathbf{E}_\sigma^\mathcal{T}1_{T_{B(\sigma,r)}>t}\mathbf{E}_{X_t}^\mathcal{T}T_{B(\sigma,r)},$$
and so substituting the above bounds for the expected exit times of balls yields\begin{equation}\label{tailstuff}
\mathbf{P}_\sigma^\mathcal{T}\left(T_{B(\sigma,r)}>t\right)\geq \frac{c_{29}}{c_{28}}\left(\ln\ln r^{-1}\right)^{-5} -\frac{t}{c_{28} r^3 \ln\ln r^{-1}}.
\end{equation}
A simple Cauchy-Schwarz argument, see \cite{Kumagai}, Proposition 4.3 for example, gives that
$$\mu(B(\sigma,r))p_{2t}(\sigma,\sigma)\geq \mathbf{P}_\sigma^\mathcal{T}\left(T_{B(\sigma,r)}>t\right)^2.$$
Now, by the volume asymptotics of Theorem \ref{localstate}, if we choose $r_{3}$ small enough, there exists a constant $c_{30}$ such that $\mu(B(\sigma,r))\leq c_{30}r^2\ln\ln r^{-1}$, for $r\in(0,r_3)$. Set $t_3= \frac{c_{29}}{2} r_{3}^{3}(\ln\ln r_{3}^{-1})^{-4}$. For $t\in(0,t_3)$, we can choose $r\in(0,r_{3})$ such that $t=\frac{c_{29}}{2} r^{3}(\ln\ln r^{-1})^{-4}$. Hence the lower bound for the tail of the exit time distribution at (\ref{tailstuff}) implies that
$$p_{2t}(\sigma,\sigma)\geq c_{31} r^{-2}\left( \ln\ln r^{-1}\right)^{-11}\geq  c_{32} t^{-2/3} \left(\ln\ln t^{-1}\right)^{-14}.$$
{\hspace*{\fill}$\square$\endlist\par\bigskip}

\hspace{10pt}\\
{\bf Proof of Theorem \ref{globalheatkernelstate}:} The upper bound of (\ref{globalheatkernelsup}) and the lower bound of (\ref{globalheatkernelinf}) are contained in Theorem \ref{process}. The lower bound of (\ref{globalheatkernelsup}) is a simple consequence of the local lower bound on the heat kernel of Theorem \ref{localheatkernelstate}. The remaining inequality is proved using the analytic technique discussed at (\ref{heatkernelupper}); the volume bound we need to utilise in this case being the lower bound for $\sup\mu(B(\sigma,r))$ appearing in Theorem \ref{globalstate}.
{\hspace*{\fill}$\square$\endlist\par\bigskip}

\section{Annealed heat kernel bounds}\label{annealedheatkernelsection}
\setcounter{thm}{0}

Rather than the $\mathbf{P}$-a.s. results about $\mu(B(\rho,r))$ and $R(\{\rho\}, B(\rho,r)^c)$ we used to prove the quenched heat kernel bounds, we need to apply estimates on the tails of their distributions to obtain annealed heat kernel bounds. We have already proved one of the necessary bounds in Lemma \ref{tailbounds0}; the remaining two bounds we require are proved in the following lemma. To complete the proof of Theorem \ref{annealedheatkernelstate}, we employ a similar argument to the proof of \cite{BarKum}, Theorem 1.4.

{\lem \label{tailbounds} There exist constants $c_{33},\dots,c_{36}$ such that, for all $r>0,\lambda\geq 1$,
$$\mathbf{P}\left(R(\{\rho\},B(\rho, r)^c)\leq r\lambda^{-1}\right)\leq c_{33}e^{-c_{34}\lambda},$$
and when $r^2\lambda^{-1}\leq \frac{1}{4}$,
$$\mathbf{P}\left(\mu(B(\rho, r))< r^2\lambda^{-1}\right)\leq c_{35}e^{-c_{36}\lambda}.$$}
\begin{proof} Let $r>0,\lambda\geq 1$. In the proof of Lemma \ref{resistance1} it was noted that $R(\{\rho\},B(\rho, r)^c)^{-1}$ is bounded above by $8r^{-1}N_{r/4}^r$. Thus, by Proposition \ref{upcrossingtail},
$$\mathbf{P}\left(R(\{\rho\},B(\rho, r)^c)\leq r\lambda^{-1}\right)\leq\mathbf{P}\left(8N_{r/4}^r\geq \lambda\right)\leq c_{16}e^{-\frac{c_{17}\lambda}{8}},$$
which proves the first inequality.

For the second inequality, suppose $r^2\lambda^{-1}\in (0,\frac{1}{4}]$. Note that if $\mu(B(\rho,r))$ is strictly less than $r^2\lambda^{-1}$, then the normalised Brownian excursion must hit the level $r$ before time $r^2\lambda^{-1}$. Thus
$$\mathbf{P}\left(\mu(B(\rho, r))< r^2\lambda^{-1}\right)\leq\mathbf{P}\left(\sup_{0\leq t\leq r^2\lambda^{-1}} W_t \geq r\right).$$
The explicit distribution of the maximum of the Brownian excursion up to a fixed time is known and we can use the formula given in \cite{DurIgl}, Section 3, to show that the right hand side of this inequality is equal to
$$1-\sqrt{\frac{2\lambda^3}{\pi r^6(1-r^2\lambda^{-1})^3}}\sum_{m=-\infty}^{\infty}e^{-2m^2r^2}\int_0^r y(2mr+y)e^{-\frac{(y+2mr(1-r^2\lambda^{-1}))^2}{2r^2\lambda^{-1}(1-r^2\lambda^{-1})}}dy.$$
We can neglect the terms with $m>0$ as removing them only increases this expression. By changing variables in the integral, it is possible to show that the $m=0$ term is equal to
$$\sqrt{\frac{2}{\pi}}\int_{0}^{\sqrt{\frac{\lambda}{1-r^2\lambda^{-1}}}}u^2 e^{-\frac{u^2}{2}}du.$$
Integrating by parts and applying standard bounds for the error function, it is elementary to obtain that this is bounded below by $1-c_{37}\sqrt{\lambda}e^{-\frac{\lambda}{2}}$. For the remaining terms we have
\begin{eqnarray*}
\lefteqn{-\sqrt{\frac{2\lambda^3}{\pi r^6(1-r^2\lambda^{-1})^3}}\sum_{m=-\infty}^{-1}e^{-2m^2r^2}\int_0^r y(2mr+y)e^{-\frac{(y+2mr(1-r^2\lambda^{-1}))^2}{2r^2\lambda^{-1}(1-r^2\lambda^{-1})}}dy}\\
&\leq & \sqrt{\frac{2\lambda^3}{\pi r^6(1-r^2\lambda^{-1})^3}}\sum_{m=1}^{\infty}\int_0^r 2mr^2e^{-\frac{(y-2mr(1-r^2\lambda^{-1}))^2}{2r^2\lambda^{-1}(1-r^2\lambda^{-1})}}dy\hspace{60pt}\\
&\leq & c_{38} \lambda^{3/2}\sum_{m=1}^{\infty} m e^{-\frac{\lambda(3m-2)^2}{8}}.
\end{eqnarray*}
The sum in this expression may be bounded above by $c_{39}e^{-\frac{\lambda}{8}}$. Thus
$$\mathbf{P}\left(\mu(B(\rho, r))< r^2\lambda^{-1}\right)\leq c_{37}\sqrt{\lambda}e^{-\frac{\lambda}{2}} + c_{38} c_{39}\lambda^{3/2}e^{-\frac{\lambda}{8}} \leq c_{40}e^{-c_{41}\lambda},$$
for suitable choice of $c_{40}, c_{41}$.
\end{proof}

\hspace{10pt}\\
{\bf Proof of Theorem \ref{annealedheatkernelstate}:} Let $t\in(0,1)$. For $\lambda\geq 2$, define $r$ by $t=2\lambda r^3$ and $A_{\lambda,t}:=\{r^2\lambda^{-1}\leq  \mu (B(\rho,r))\leq r^2\lambda\}$. On $A_{\lambda,t}$, we can use the inequality at (\ref{heatkernelupper}) to show that
$$p_t(\rho,\rho)\leq \frac{2^{5/3}\lambda^{5/3}}{t^{2/3}}.$$
Define $\Lambda_t:=\inf\{\lambda\geq 2:\:A_{\lambda,t}\mbox{ occurs}\}$, then $\mathbf{E}p_t(\rho,\rho)\leq 2^{5/3} t^{-2/3} \mathbf{E}\Lambda_t^{5/3}$. However, for $\lambda\geq 2$,
$$\mathbf{P}(\Lambda_t\geq\lambda)\leq\mathbf{P}(A_{\lambda,t}^c)\leq\mathbf{P}(\mu(B(\rho,r))> r^2\lambda)+\mathbf{P}(\mu(B(\rho,r))< r^2\lambda^{-1}).$$
Since $r^2\lambda^{-1}=t^{2/3} 2^{-2/3} \lambda^{-5/3}\leq \frac{1}{4}$, we can apply the tail bounds of Lemmas \ref{tailbounds0} and \ref{tailbounds} to obtain that $\mathbf{P}(\Lambda_t\geq \lambda)\leq c_{42}e^{-c_{43}\lambda}$, uniformly in $t\in (0,1)$. Thus $\mathbf{E}\Lambda_t^{5/3}\leq c_{44}<\infty$, uniformly in $t\in(0,1)$, which proves the upper bound.

For the lower bound we need a slightly different scaling. Let $t\in(0,1)$, $\lambda\geq 64$ and define $r$ by $t=r^3/32\lambda^4$, and
$$B_{\lambda,t}:=\left\{\mu(B(\rho,r))\leq r^2\lambda,\:R(\{\rho\},B(\rho,r)^c)\geq r\lambda^{-1},\:\mu(B(\rho,\frac{r}{4\lambda}))\geq \frac{r^2}{16\lambda^3}\right\}.$$
On $B_{\lambda,t}$, by following a similar argument to that used for the proof of Theorem \ref{localheatkernelstate}, we find that $p_t(\rho,\rho)\geq c_{45} t^{-2/3}\lambda^{-14}$. Now,
\begin{eqnarray*}
\mathbf{P}(B_{\lambda,t}^c)&\leq &\mathbf{P}(\mu(B(\rho,r))> r^2\lambda)+\mathbf{P}(R(\{\rho\},B(\rho,r)^c)< r\lambda^{-1})\\
& & \hspace{140pt}+\mathbf{P}(\mu(B(\rho,\frac{r}{4\lambda}))< \frac{r^2}{16\lambda^3}).
\end{eqnarray*}
Since $r^2/16\lambda^{3}=t^{2/3}\lambda^{-1/3}\leq \frac{1}{4}$, again we can apply the bounds of Lemmas \ref{tailbounds0} and \ref{tailbounds} to find that $\mathbf{P}(B_{\lambda,t}^c)\leq c_{46}e^{-c_{47}\lambda}$, uniformly in $t\in(0,1)$. Hence we can find a $\lambda_0\in[64,\infty)$ such that $\mathbf{P}(B_{\lambda_0,t}^c)\leq\frac{1}{2}$ for all $t\in(0,1)$. Thus
$$\mathbf{E}p_t(\rho,\rho)\geq \mathbf{P}(B_{\lambda_0,t})\frac{c_{45}}{\lambda_0^{14}t^{2/3}}\geq c_{48}t^{-2/3},\hspace{20pt}\forall t\in(0,1),$$
for some $c_{48}>0$.
{\hspace*{\fill}$\square$\endlist\par\bigskip}

\section{Brownian motion on the CRT}\label{BM}
\setcounter{thm}{0}

To complete the proof of Theorem \ref{heatkernel}, it remains to show that the Markov process with infinitesimal generator $\Delta_\mathcal{T}$ is Brownian motion on $\mathcal{T}$. Brownian motion on $\mathcal{T}_f$ is defined to be a $\mathcal{T}_f$-valued process, $X^f=((X^f_t)_{t\geq 0}, \mathbf{P}^{\mathcal{T}_f}_\sigma, \sigma\in\mathcal{T}_f)$, with the following properties.
\newcounter{listcount}
\begin{list}{\roman{listcount})}{\usecounter{listcount} \setlength{\rightmargin}{\leftmargin}}
\item Continuous sample paths.
\item Strong Markov.
\item Reversible with respect to its invariant measure $\mu_f$.
\item For $\sigma^1,\sigma^2\in\mathcal{T}_f$, $\sigma^1\neq\sigma^2$, we have
$$\mathbf{P}_{\sigma}^{\mathcal{T}_f}\left( T_{\sigma^1}<T_{\sigma^2}\right)=\frac{d_{\mathcal{T}_f}(b(\sigma,\sigma^1,\sigma^2),\sigma^2)}{d_{\mathcal{T}_f}(\sigma^1,\sigma^2)},\hspace{20pt}\sigma\in\mathcal{T}_f,$$
where $T_\sigma:=\inf\{t\geq 0:\:X^f_t=\sigma\}$ and $b(\sigma,\sigma^1,\sigma^2)$ is the branch point of $\sigma, \sigma^1,\sigma^2$ in $\mathcal{T}$. In particular, if $[[\sigma,\sigma^1]]$, $[[\sigma^1,\sigma^2]]$,  $[[\sigma^2,\sigma]]$ are the unique arcs between the relevant pairs of vertices, then $b(\sigma,\sigma^1,\sigma^2)$ is the unique point in the set $[[\sigma,\sigma^1]]\cap[[\sigma^1,\sigma^2]]\cap[[\sigma^2,\sigma]]$.
\item For $\sigma^1,\sigma^2\in\mathcal{T}_f$, the mean occupation measure for the process started at $\sigma^1$ and killed on hitting $\sigma^2$ has density
$$d_{\mathcal{T}_f}(b(\sigma,\sigma^1,\sigma^2),\sigma^2)\mu(d\sigma),\hspace{20pt}\sigma\in\mathcal{T}_f.$$
\end{list}
As remarked in \cite{Aldous2}, Section 5.2, these properties are enough for uniqueness of Brownian motion on $\mathcal{T}_f$. Note that the definition given by Aldous has an extra factor of 2 in property v). This is a result of Aldous' description of the CRT being based on the random function $2W$. By Theorem \ref{process}, we already have that properties i), ii) and iii) hold for the process $X$ on $\mathcal{T}$, $\mathbf{P}$-a.s.

An important tool for the proofs of properties iv) and v) will be the trace operator for Dirichlet forms, which will allow us to deduce results about the Dirichlet form on $\mathcal{T}$ by considering its restriction to finite subsets of $\mathcal{T}$. We now introduce the notation for this, for more details, see \cite{Barlow}, Section 4. Let $K$ be a set and $\nu$ a measure on $K$ and suppose we are given a regular Dirichlet form $(\mathcal{E},\mathcal{F})$ on $L^2(K,\nu)$. Let $\tilde{K}$ be the closed support of the $\sigma$-finite measure $\tilde{\nu}$ on $K$. The trace of $\mathcal{E}$ on $\tilde{K}$ is defined by
$${\rm Tr}(\mathcal{E}|\tilde{K})(v,v):=\inf\{\mathcal{E}(u,u):\:u|_{\tilde{K}}=v\},\hspace{20pt}v\in L^2(\tilde{K},\tilde{\nu}),$$
and its domain, $\tilde{\mathcal{F}}$, is the set of functions for which this infimum is finite. Before proceeding with the main results of this section, we first need to prove the following technical result on the capacity of sets of $\mathcal{T}$, where we use the notation $\tilde{U}$ to represent the set of excursions on which a finite resistance form is defined on $\mathcal{T}_f$, as in Section \ref{hkbounds}.

{\lem For $f\in \tilde{U}$, all non-empty subsets of $\mathcal{T}_f$ have strictly positive capacity.}
\begin{proof}
For $f\in\tilde{U}$, we are able to construct the associated Dirichlet form, $\mathcal{E}=\mathcal{E}_{\mathcal{T}_f}$. By Lemma 3.2.2 of \cite{Fukushima}, if $\nu$ is a positive Radon measure on $\mathcal{T}_f$ with finite energy integral, i.e.,
$$\left(\int_{\mathcal{T}_f}|u(\sigma)|\nu(d\sigma)\right)^2\leq c_{49} \left(\mathcal{E}(u,u)+\int_{\mathcal{T}_f} u(\sigma)^2\mu(d\sigma)\right),\hspace{20pt}\forall u\in \mathcal{F},$$
for some $c_{49}<\infty$, then $\nu$ charges no set of zero capacity. Hence, by monotonicity, it is sufficient to show that $\delta_\sigma$ has finite energy integral for any $\sigma\in\mathcal{T}_f$, where $\delta_{\sigma}$ is the probability measure putting all its mass at the point $\sigma$. For any $\sigma,\sigma'\in\mathcal{T}_f$,
$$\left(\int_{\mathcal{T}_f}|u(\sigma'')|\delta_\sigma(d\sigma'')\right)^2= u(\sigma)^2\leq  2 (u(\sigma)-u(\sigma'))^2+2u(\sigma')^2.$$
Applying the inequality at (\ref{contin}) to this bound, and integrating with respect to $\sigma'$ yields
$$\left(\int_{\mathcal{T}_f}|u(\sigma')|\delta_\sigma(d\sigma')\right)^2\leq 2\:{\rm diam}\mathcal{T}_f\ \mathcal{E}(u,u)+2\int_{\mathcal{T}_f} u(\sigma')^2\mu(d\sigma'),$$
which completes the proof, because ${\rm diam}\mathcal{T}_f$ is finite.
\end{proof}

The above result allows us to define local times for $X$, and in the following lemma we use these to define a time-changed process on a finite subset of $\mathcal{T}$. By considering the hitting probabilities for the time-changed process, we are able to deduce that $X$ satisfies property iv) of the Brownian motion definition. In particular, we will use the fact that the quadratic form corresponding to our time-changed process is simply the trace of $\mathcal{E}_\mathcal{T}$ on the same finite subset.

{\lem \label{prop4} $\mathbf{P}$-a.s., the process $X$ of Theorem \ref{process} satisfies property iv) of the definition of Brownian motion on $\mathcal{T}$.}
\begin{proof}
Suppose $W\in \tilde{U}$, so that the resistance form, $\mathcal{E}_\mathcal{T}$, and process, $X$, are defined for $\mathcal{T}$. Fix $\sigma, \sigma^1,\sigma^2\in\mathcal{T}$, $\sigma^1\neq\sigma^2$, and set $b=b(\sigma,\sigma^1,\sigma^2)$. Write $V_1=\{\sigma,\sigma^1,\sigma^2,b\}$ and $\mathcal{E}_1={\rm Tr}(\mathcal{E}_\mathcal{T}|V_1)$. Using simple properties of resistance forms, the following explicit expression for $\mathcal{E}_1$ can be calculated:
\begin{equation}\label{resform1}
\mathcal{E}_1(u,u)=\sum_{\sigma'\in\{\sigma,\sigma^1,\sigma^2\}} \frac{\left(u(b)-u(\sigma')\right)^2}{d_\mathcal{T}(b, \sigma')},\hspace{20pt}u\in C(V_1),
\end{equation}
where, if $b=\sigma'$, the relevant term is defined to be 0.

By the previous lemma, $\{\sigma'\}$ has strictly positive capacity for each $\sigma'\in\mathcal{T}$. As outlined in Section 4 of \cite{Barlow}, a result of this is that $X$ has jointly measurable local times $(L_t^{\sigma'},\sigma'\in\mathcal{T},t\geq 0)$ such that
$$\int_0^tu(X_s)ds=\int_\mathcal{T}u(\sigma')L_t^{\sigma'}\mu(d\sigma'),\hspace{20pt}u\in L^2(\mathcal{T},\mu).$$
Now, denote $\nu:=\frac{1}{|V_1|}\sum_{\sigma'\in V_1} \delta_{\sigma'}$, the uniform distribution on $V_1$ and define
$$A_t:=\int_\mathcal{T} L_t^{\sigma'} \nu(d\sigma'),\hspace{20pt}\tau_t:=\inf\{s:\:A_s>t\}.$$
Consider the process $\tilde{X}=(\tilde{X},\mathbf{P}_{\sigma'}^\mathcal{T},\sigma'\in V_1)$, defined by $\tilde{X}_t:=X_{\tau_t}$. As described in \cite{Barlow}, $\tilde{X}$ is a $\nu$-symmetric Hunt process and has associated regular Dirichlet form $(\mathcal{E}_1, C(V_1))$. Using elementary theory for continuous time Markov chains on a finite state space, we obtain the following result for $\tilde{X}$
$$\mathbf{P}_\sigma^{\mathcal{T}} \left(\tilde{T}_{\sigma^1}<\tilde{T}_{\sigma^2}\right)= \frac{d_{\mathcal{T}}(b,\sigma^2)}{d_{\mathcal{T}}(\sigma^1,\sigma^2)},$$
where $\tilde{T}_{\sigma'}:=\inf\{t\geq 0:\:\tilde{X}_t=\sigma'\}$. Since the hitting distribution is unaffected by the time change from $X$ to $\tilde{X}$, this implies
$$\mathbf{P}_\sigma^{\mathcal{T}} \left({T}_{\sigma^1}<{T}_{\sigma^2}\right)= \frac{d_{\mathcal{T}}(b,\sigma^2)}{d_{\mathcal{T}}(\sigma^1,\sigma^2)},$$
and so, if $W\in \tilde{U}$, the process $X$ satisfies property iv) of the Brownian motion definition. Since $W\in \tilde{U}$, $\mathbf{P}$-a.s., this completes the proof.
\end{proof}

A result that will be useful in proving that $X$ satisfies property v) is the following uniqueness result, which is proved in \cite{Kigamidendrite}, Lemma 3.5.

{\lem \label{unique} Let $(\mathcal{E},\mathcal{F})$ be a resistance form on a set $K$ and $V$ be a finite subset of $K$. Then for any $v\in C(V)$, there exists a unique $u\in \mathcal{F}$ such that
$$\mathcal{E}(u,u)={\rm Tr}(\mathcal{E}|V)(v,v),\hspace{20pt}u|_V=v.$$}

{\lem \label{prop5} $\mathbf{P}$-a.s., the process $X$ of Theorem \ref{process} satisfies property v) of the definition of Brownian motion on $\mathcal{T}$.}
\begin{proof}
First, assume $W\in \tilde{U}$, so that the resistance form, $\mathcal{E}_\mathcal{T}$, and process, $X$, are defined for $\mathcal{T}$. Fix $\sigma^1,\sigma^2\in\mathcal{T}$, $\sigma^1\neq\sigma^2$, and define $D=D(\sigma^1,\sigma^2)$ to be the path connected component of $\mathcal{T}\backslash\{\sigma^2\}$ containing $\sigma^1$. Using the same argument as in \cite{Kumagai}, Proposition 4.2, we can deduce the existence of a Green kernel $g^D(\cdot, \cdot)$ for the process killed on exiting $D$ which satisfies
\begin{equation}\label{reproducing}
\mathcal{E}_\mathcal{T}(g^D(\sigma,\cdot),f)=f(\sigma),\hspace{20pt}\forall \sigma\in\mathcal{T}, f\in \mathcal{F}_D,
\end{equation}
where $\mathcal{F}_D:=\{f\in\mathcal{F}_\mathcal{T}:\:f|_{D^c}=0\}$. By standard arguments, this implies that $g^D(\sigma^1,\sigma^1)>0$;
\begin{equation}\label{infresult}
\mathcal{E}_\mathcal{T}(\tilde{g},\tilde{g})=\inf\{\mathcal{E}_\mathcal{T}(u,u):\:u(\sigma^1)=1,\:u(\sigma^2)=0\},
\end{equation}
where $\tilde{g}(\cdot):=g^D(\sigma^1,\cdot)/g^D(\sigma^1,\sigma^1)$; and for $\mu$-measurable $f$,
$$\mathbf{E}^{\sigma^1}\int_0^{T_{\sigma^2}}f(X_s)ds=\int_\mathcal{T} g^D(\sigma^1,\sigma)f(\sigma)\mu(d\sigma).$$
This means that $g^D(\sigma^1,\sigma)\mu(d\sigma)$ is the mean occupation density of the process started at $\sigma^1$ and killed on hitting $\sigma^2$. Furthermore, note that combining (\ref{reproducing}), (\ref{infresult}) and the characterisation of $d_\mathcal{T}$ at (\ref{resistance}), we can deduce that $g^D(\sigma^1,\sigma^1)=d_\mathcal{T}(\sigma^1,\sigma^2)$.

Now, fix $\sigma\in\mathcal{T}$, and define $b:=b(\sigma,\sigma^1,\sigma^2)$,
$$V_1:=\{\sigma,\sigma^1,\sigma^2,b \},\hspace{20pt}\mathcal{E}_1:={\rm Tr}(\mathcal{E}_\mathcal{T}|V_1),$$
$$V_0:=\{\sigma^1,\sigma^2 \},\hspace{20pt}\mathcal{E}_0:={\rm Tr}(\mathcal{E}_1|V_0).$$
Let $f_0\in C(V_0)$ be defined by $f_0(\sigma^1)=1$, $f_0(\sigma^2)=0$; $f_1$ be the unique (by Lemma \ref{unique}) function in $C(V_1)$ that satisfies $\mathcal{E}_1(f_1,f_1)=\mathcal{E}_0(f_0,f_0)$ and $f_1|_{V_0}=f_0$; and $f_2$ be the unique function in $\mathcal{F}$ such that $\mathcal{E}_\mathcal{T}(f_2,f_2)=\mathcal{E}_1(f_1,f_1)$ and $f_2|_{V_1}=f_1$. Applying the tower property for the trace operator, $\mathcal{E}_0 ={\rm Tr}({\rm Tr}(\mathcal{E}_\mathcal{T}|V_1)|V_0)={\rm Tr}(\mathcal{E}_\mathcal{T}|V_0)$, we have that $f_2$ is the unique function that satisfies
$$\mathcal{E}_\mathcal{T}(f_2,f_2)={\rm Tr}(\mathcal{E}_\mathcal{T}|V_0)(f_0,f_0),\hspace{20pt}f_2|_{V_0}=f_0.$$
However, we have from (\ref{infresult}) that $\tilde{g}$ also has these properties and so it follows from the uniqueness of Lemma \ref{unique} that $\tilde{g}=f_2$. Thus $\tilde{g}|_{V_1}=f_1$. Recall, the explicit expression for $\mathcal{E}_1$ given at (\ref{resform1}). A simple minimisation of this quadratic polynomial allows us to determine the function $f_1$. In particular, we have $\tilde{g}(\sigma)=f_1(\sigma)=d_{\mathcal{T}}(b,\sigma^2)d_\mathcal{T}(\sigma^1,\sigma^2)^{-1}$. Hence the mean occupation density of the process started at $\sigma^1$ and killed on hitting $\sigma^2$ is
$$g^D(\sigma^1,\sigma)\mu(d\sigma)= g^D(\sigma^1,\sigma^1)\tilde{g}(\sigma)\mu(d\sigma)= d_{\mathcal{T}}(b,\sigma^2)\mu(d\sigma).$$
Thus, if $W\in \tilde{U}$, the process $X$ satisfies property v) of the Brownian motion definition. Since $W\in \tilde{U}$, $\mathbf{P}$-a.s., this completes the proof.
\end{proof}

Combining the results of Theorem \ref{process} and Lemmas \ref{prop4} and \ref{prop5} we immediately have the following.

{\cor $\mathbf{P}$-a.s., the process $X$ of Theorem \ref{process} is Brownian motion on $\mathcal{T}$.}

\def\cprime{$'$}
\providecommand{\bysame}{\leavevmode\hbox to3em{\hrulefill}\thinspace}
\providecommand{\MR}{\relax\ifhmode\unskip\space\fi MR }
\providecommand{\MRhref}[2]{%
  \href{http://www.ams.org/mathscinet-getitem?mr=#1}{#2}
}
\providecommand{\href}[2]{#2}

\end{document}